\documentclass[12pt]{amsart}

\usepackage{amsfonts,amssymb,amscd,epsf,xy}

\xyoption{all}

\def\dual                 {{\vee}}

\def\ee                 {{\rm e}}

\def\Box             {{\rm Box}}
\def\sup            {{\rm Supp}}

\def\ZZ                 {{\mathbb Z}}
\def\PP                {{\mathbb P}}
\def\RR                 {{\mathbb R}}
\def\CC                 {{\mathbb C}}
\def\QQ                 {{\mathbb Q}}
\def\LL                 {{\mathbb L}}

\def\Aa    {{\mathcal A}}
\def\Ba    {{\mathcal B}}
\def\Ca    {{\mathcal C}}
\def\Da    {{\mathcal D}}
\def\Fa    {{\mathcal F}}
\def\Ia    {{\mathcal I}}

\def\Ra    {{\mathcal R}}
\def\Sa    {{\mathcal S}}
\def\Ta    {{\mathcal T}}

\def\Rh    {{\hat{R}}}
\def\Rah   {{\hat{\mathcal R}}}

\def\th    {{\hat{t}}}

\newtheorem{lemma}{Lemma}[section]
\newtheorem{theorem}[lemma]{Theorem}
\newtheorem{corollary}[lemma]{Corollary}
\newtheorem{proposition}[lemma]{Proposition}
\theoremstyle{definition}
\newtheorem{definition}[lemma]{Definition}

\newtheorem{remark}[lemma]{Remark}
\theoremstyle{remark}
\newtheorem*{proof*}{Proof}

\title{Mellin--Barnes integrals as Fourier--Mukai transforms}

\author{Lev A. Borisov and R. Paul Horja}
\address{Department of Mathematics \\ University of Wisconsin \\
  Madison \\ WI \\ 53706 \\ USA\\{\tt borisov@math.wisc.edu}}
\address{Department of Mathematics \\ Oklahoma State University \\
  Stillwater \\ OK \\ 74078 \\ USA\\{\tt horja@math.okstate.edu}}

\thanks{The first 
author was partially supported by NSF grant DMS-0140172.}

\begin{document}

\begin{abstract}
We study the generalized hypergeometric 
system introduced by Gelfand, Kapranov and Zelevinsky 
and its relationship with the toric Deligne--Mumford (DM)
stacks recently studied by Borisov, Chen and Smith. We construct
series solutions with values in a combinatorial version of the
Chen--Ruan (orbifold)
cohomology and in the $K$--theory of the associated DM stacks. In the spirit 
of the homological
mirror symmetry conjecture of Kontsevich, we show that the $K$--theory 
action of the Fourier--Mukai functors associated to basic
toric birational maps of DM stacks 
are mirrored by analytic continuation transformations of 
Mellin--Barnes type. 
\end{abstract}

\maketitle

\section{Introduction}

Let $\Aa= \{ v_1, \ldots, v_n \}$ be a collection of elements of the 
lattice $N \cong \ZZ^d.$ We assume that the elements of $\Aa$ 
generate the lattice as an abelian group, and that there exists a
group homomorphism $h : N \to \ZZ$ such that $h(v)=1$ for any element
$v \in \Aa.$ Let $\LL \subset \ZZ^n$ denote the lattice of integral relations 
among the elements of $\Aa$ consisting of vectors $l=(l_j) \in \ZZ^n$ 
such that $l_1 v_1 + \ldots + l_n v_n =0.$ 

Let $\beta \in N$ be a lattice element. 
The {\it Gelfand--Kapranov--Zelevinsky hypergeometric system (GKZ)} associated
to the set $\Aa$ and parameter $\beta$ 
is a 
system of differential equations on the function $\Phi(z),$
$z=(z_1, \ldots, z_n) \in \CC^n,$ consisting of the
binomial equations
$$
\Big( \prod_{j, l_j >0} \big(\frac{\partial}{\partial z_j} \big)^{l_j}
- \prod_{j, l_j < 0} \big(\frac{\partial}{\partial z_j} \big)^{-l_j} \Big)
\Phi =0, \ l \in \LL, 
$$
and the linear equations
$$
\Big(-\beta + \sum_{j=1}^n v_j z_j \frac{\partial}{\partial z_j} \Big) \Phi= 0.
$$
Note that it is enough to consider a finite set of binomial equations determined
by a set of generators of the lattice of relations $\LL.$ 

Gelfand, Kapranov and Zelevinsky \cite{GKZ1}
showed that this system is holonomic, so the number
of solutions at a generic point is finite. They constructed explicit 
solutions of the system in the form of the so-called Gamma series
$$
\Phi(z, \lambda):= \sum_{l \in \LL} \prod_{j=1}^n \frac{z_j^{l_j + \lambda_j}}
{\Gamma(l_j + \lambda_j +1)},
$$
where $\lambda \in \CC^n$ is a parameter with the property that 
$\lambda_1 v_1 + \ldots + \lambda_n v_n = \beta.$ Moreover, 
they discovered that there is a very close connection
between the regular triangulations 
of the polytope $\Delta={\rm Conv}(\Aa),$ 
as described by the secondary polytope
of $\Aa,$ and the structure of the solution set. 

In the context of mirror symmetry, 
Batyrev \cite{Bat} noticed that a special case of the 
GKZ system is satisfied by the periods 
describing the variations of complex structures of  
Calabi--Yau hypersurfaces in toric varieties. Aspinwall,
Greene and Morrison \cite{AGM} used the combinatorial GKZ machinery 
to analyze the string theoretic phase transitions
in type II string theory. 
The homological
mirror symmetry conjecture of Kontsevich \cite{Kont1}
provided a far reaching generalization of the earlier
ideas in mirror symmetry. As further evidence for his 
proposal, Kontsevich \cite{Kont2} conjectured that the 
action on cohomology of the group of 
self--equivalences of the bounded derived category 
of coherent sheaves on a smooth projective Calabi--Yau
variety matches the monodromy action on the cohomology
of the mirror Calabi--Yau variety associated to the 
variations of complex structures. In the toric context,
this strategy has been pursued by one of the authors 
in \cite{H}. The broad goal 
of the current work is to 
offer a framework for the aforementioned ideas based on the
notion of a toric Deligne--Mumford stack introduced
by Borisov, Chen and Smith \cite{BCS}. 

Sections \ref{sec:gkzvalsr} and \ref{sec:gkzk} 
provide a geometric 
approach to the problem of
constructing convergent Gamma series solutions to the GKZ system 
corresponding to a general regular triangulation of the 
polytope $\Delta={\rm Conv}(\Aa).$ 
For an extensive list of works that investigate
the properties of the GKZ system and the associated $\Da$--module,
the bibliography of the book by Saito, Sturmfels and Takayama 
\cite{SST} is the best resource. 
Our approach is closest in spirit to the methods employed by 
Hosono, Lian and Yau \cite{HLY} and 
Stienstra \cite{stienstra} for the case of unimodular
triangulations which, for some time, has been the preferred 
testing ground
for mirror symmetry computations. However, 
homological
mirror symmetry and the advent of D-branes in string theory have sparked 
renewed interest in understanding the intricacies of the general
situation. What we show is that the  
Chen-Ruan (orbifold) cohomology and the $K$--theory of the toric
Deligne--Mumford stack associated to a general regular triangulation 
of the the convex polytope $\Delta$
are the natural missing ingredients in the 
geometric construction of the solution set to the GKZ system. 

In particular, for a fan $\Sigma$ supported on the cone
$K=\RR_{\geq 0} \Delta$
induced by a regular triangulation 
of the polytope $\Delta,$ and an arbitrary $\beta \in N,$ 
we explicitly construct at least ${\rm Vol}(\Delta)$
cohomology valued linearly independent Gamma series solutions 
in a certain complex domain 
$U_{\Sigma}$ in $\CC^n$ associated to the fan $\Sigma$
(corollary \ref{cor:gkzsvol}$i)$). 
When $\beta \in -K^{\circ},$ 
which is the case of mirror symmetry 
for Calabi-Yau complete intersections in projective
toric varieties, the second part of the same corollary provides 
a full system of ${\rm Vol}(\Delta)$ linearly independent 
solutions. 
Furthermore, corollary \ref{cor:msmap} provides a mirror
symmetry map defined on the dual of the Grothendieck ring of the 
toric DM stack $\PP_\Sigma$ 
$$
MS_{\Sigma} : (K_0(\PP_{\Sigma}, \CC))^{\dual} 
\to {\mathcal Sol} (U_{\Sigma})
$$
that produces 
GKZ solutions which are analytic in the complex domain 
$U_{\Sigma}$ in $\CC^n$ associated to the fan $\Sigma.$

Our $K$--theoretic interpretation
makes essential use of results contained in the companion paper \cite{BH},
where, among other things, a Stanley--Reisner type description of the 
Grothendieck $K$--theory ring of a smooth DM stack
is given. It is well known that there are many subtleties involved 
in trying to determine 
the dimension of the solution set of the GKZ system (see, for example, 
Adolphson \cite{A}, Saito, Sturmfels and Takayama \cite{SST}, 
Cattani, Dickenstein and Sturmfels \cite{CDS}, Matusevich, Miller 
and Walther \cite{MMW}). 
Our methods raise the interesting issue 
of finding the proper $K$--theoretic framework for constructing 
GKZ solutions for general values of the parameter $\beta \in N.$

In sections \ref{sec:achs} and \ref{sec:vs}, we employ 
a combination of analytic, algebra-geometric and combinatorial methods 
and justify the title of the paper. We consider $\Sigma_+$
and $\Sigma_-$ two fan structures induced by two regular
triangulations of the polytope $\Delta$ that are joined by an 
edge of the secondary polytope determined by $\Aa.$ It follows 
that the associated toric DM stacks $\PP_{\Sigma_+}$ and
$\PP_{\Sigma_-}$ are birationally equivalent, and 
we have a diagram of weighted blowdowns 
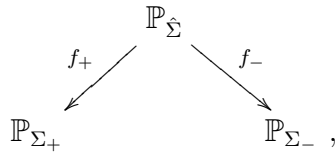
\begin{equation*}
\begin{split}
\xymatrix{
& \PP_{\hat{\Sigma}} \ar[ld]_{f_+} 
 \ar[rd]^{f_-} \\
\PP_{\Sigma_+} & & \PP_{\Sigma_-} \ , 
}
\end{split}
\end{equation*}
where $\hat{\Sigma}$ is a stacky refinement of the fans $\Sigma_\pm.$
According to 
Bondal 
and Orlov \cite{BO}, in the smooth fan case, and Kawamata \cite{Ka},
in the stacky case, the map $FM : {\rm D}^b(\PP_{\Sigma_-}) \to 
{\rm D}^b(\PP_{\Sigma_+})$ between the bounded derived categories 
of coherent sheaves of the associated DM stacks, 
defined by
$$FM:=  {\bf R} (f_+)_* \, {\bf L}(f_-)^*,$$ 
is an equivalence of triangulated categories, i.e. a Fourier--Mukai
functor. Our main result, theorem \ref{thm:comm}, shows 
the commutativity of the diagram
\begin{equation*}
\begin{split}
\xymatrix@C=1in{
(K_0(\PP_{\Sigma_+},\CC))^\vee \ar[r]^{MS_{+}}
 \ar[d]_{FM^\vee}
 & {\mathcal Sol}(U_{+})\ar[d]^{MB}\\
(K_0(\PP_{\Sigma_-},\CC))^\vee \ar[r]^{MS_{-}}
 & {\mathcal Sol}(U_{-})
}
\end{split}
\end{equation*}
Not surprisingly, in accordance with general mirror symmetry
principles, 
the left side of the diagram has to do with birational geometry, 
while its right side is analytic in nature. The explicit 
description of the analytic continuation map $MB$ is given in
corollary \ref{cor:mb} (see definition \ref{def:mb}) 
as an application of the Mellin--Barnes integral representation  method. 
The use of 
$K$--theory instead of (orbifold) cohomology in this last 
part of the paper has many advantages. 
In particular, we do not need to use any type of 
Grothendieck--Riemann--Roch theorem for stacks. 

Most of the results of this work were obtained
while the second author held
positions at the University of Michigan, Ann Arbor, 
and at the Fields Institute, University of Toronto. 
He is very grateful for financial support and for 
the outstanding scientific atmosphere provided by
both institutions. 
We would like to thank Victor Batyrev and Kentaro Hori 
for useful discussions. 

{\it Notation.} 
Given a fan $\Sigma$ in $N,$ and two cones $\sigma, \tau \in \Sigma,$ 
we write $\sigma \prec \tau$ or $\tau \succ \sigma$
to indicate that the cone $\sigma$ is a face of the cone 
$\tau.$ A one-dimensional face of a cone $\sigma$ will
sometimes be called a \emph{ray} of $\sigma.$
For a subset $\Ba$ of the lattice $N,$
we write $\RR_{\geq 0} \Ba$ to denote the cone 
generated by the elements of $\Ba.$

\section{GKZ solutions with values in toric SR cohomology}
\label{sec:gkzvalsr}

As in the previous section, assume that
$\Aa= \{ v_1, \ldots, v_n \} \subset N \cong \ZZ^d$ 
generates the lattice $N,$ and that all the elements of 
$\Aa$ are located in a hyperplane $h (w)=1,$ for 
a linear map $h: N \to \ZZ.$ 
In what follows, we consider regular
triangulations of the polytope $\Delta= {\rm Conv} (\Aa)$
with all their vertices among the elements of 
$\Aa.$ Every such triangulation 
determines a fan structure $\Sigma$ supported on the cone 
$K=\RR_{\geq 0} \Delta.$
It is well known (see chapter 7 in \cite{GKZbook}) 
that there exists a one-to-one correspondence
between the regular triangulations of $\Delta$ 
and the maximal cones of the secondary fan determined by $\Aa$
(or, dually, the vertices of the secondary polytope associated to 
$\Aa$). 

We define the partial semigroup ring $\CC[K,\Sigma]$ associated to the
cone $K$ and the fan $\Sigma$ to be the complex vector space with a basis
given by the symbols $x^w$ for all $w \in K \cap N$ and the multiplication
defined such that $x^{w_1} \cdot x^{w_2}= x^{w_1 + w_2},$ whenever there
exists a cone $\sigma \in \Sigma$ containing both $w_1$ and $w_2,$ and 
$x^{w_1} \cdot x^{w_2}= 0,$ otherwise. 
The ideal $\CC[K^\circ,\Sigma]
\subset \CC[K,\Sigma]$ associated to the interior of the cone $K$ 
is generated by the elements $x^w$ for all the elements $w \in K^\circ.$
The ring $\CC[K,\Sigma]$ 
and the ideal $\CC[K^\circ,\Sigma]$
admit a natural positive grading induced by the hyperplane condition 
on $\Aa.$ 

The following result is stated in \cite{BM} (see also \cite{B}). 

\begin{proposition}\label{prop:regseq}
The ring $\CC[K,\Sigma]$ and the module $\CC[K^\circ,\Sigma]$ 
(over  $\CC[K,\Sigma]$) are Cohen-Macaulay of dimension $d.$ 
Moreover,
for any basis $(m_1,\ldots, m_d )$ of $M={\rm Hom} (N, \ZZ),$ the elements
$$
Z_i=\sum_{j,\RR_{\geq 0}v_j \in \Sigma} \langle m_i, v_j \rangle x^{v_j}
$$
form a regular sequence in $\CC[K,\Sigma]$ 
(and hence in $\CC[K^\circ,\Sigma]$). 
\end{proposition}

\begin{corollary}
The quotients  
$$
\CC[K,\Sigma]/Z \CC[K,\Sigma]:= 
\CC[K,\Sigma]/(Z_1,\ldots,Z_d)\CC[K,\Sigma]$$ 
and 
$$
\CC[K^\circ,\Sigma]/Z \CC[K^\circ,\Sigma]:=
\CC[K^\circ,\Sigma]/(Z_1,\ldots,Z_d)\CC[K^\circ,\Sigma]$$
have dimension equal to the normalized volume of $\Delta$.
\end{corollary}
\begin{proof} 
The dimensions of these vector spaces are 
equal to $(d-1)!$ times the leading coefficient
of the Hilbert polynomial of the graded ring 
and module. It is well-known this leading
coefficient  is 
the quotient of the normalized volume by $(d-1)!,$
see for example \cite{Sturm}, theorem 4.16.
\end{proof}

In line with the terminology of section 3 in \cite{BH}, we will call 
the quotient $\CC[K,\Sigma]/Z\CC[K,\Sigma]$ 
the {\it SR--cohomology ring.} 

\subsection{$\Gamma$--series with values in the completion of
$\CC[K,\Sigma]$}
\label{subsec:cke}

For any maximal dimensional cone $\sigma \in \Sigma,$ 
we define the 
set $\Box(\sigma)$ of elements in $N$ to be the set
$$
\{ v \, :  \, v = \sum_{j=1}^n
q^{v}_j v_j, \ 0 \leq q^{v}_j <1, \ q^v_j=0, \ 
\text{if} \ \RR_{\geq 0} v_j \ \text{is not a ray of} \ \sigma  \}.
$$
Define the set $\Box (\Sigma)$ of elements in $N$ to be 
the union of the sets $\Box(\sigma)$ for all the maximal 
dimensional cones $\sigma \in \Sigma.$ If $v \in 
\Box(\Sigma),$ we denote by $\sigma(v)$ 
the smallest cone of $\Sigma$ that contains $v.$ 

The set $\mathcal A$ generates the lattice $N,$
so for each $v \in \Box(\Sigma),$ we can choose a 
solution $\gamma^v$
of the equation $\gamma^{v}_1 v_1+\ldots+\gamma^{v}_n v_n
=\beta,$
with the property that $\gamma^{v}_j \equiv q^{v}_j$ (mod 
$\mathbb Z$). This implies that $\gamma_j^v$ is integer, unless
$\RR_{\geq 0} v_j$ is a ray of $\sigma(v).$
In particular, 
$\gamma_j^v$ is integer for $\RR_{\geq 0} v_j\not\in \Sigma$.

For a given $z=(z_1, \ldots, z_n)$ in $(\CC^\star)^n,$
consider the formal expression 
$$
x^{v} \prod_{j=1}^n  \frac{z_j^{l_j + \gamma^v_j+D_j}}
{\Gamma(l_j+ \gamma^v_j+D_j+1)},
$$
where $l \in \LL,$
$$
D_j := x^{v_j} \ \text{if} \ \RR_{\geq 0} v_j \in \Sigma, \ \text{and} \
D_j:=0, \ \text{otherwise,}$$ 
and 
$$
z_j^{\gamma^v_j+D_j}:= e^{(\gamma^v_j+D_j) (\log |z_j| + i \arg z_j)},$$
for a choice of $(\arg z_1, \ldots, \arg z_n) \in \RR^n.$

\begin{definition}\label{def:coner}
Let $\gamma=(\gamma_1, \ldots, \gamma_n)$ be an 
element in $\QQ^n.$ For any $l=(l_1, \ldots, l_n) \in \LL,$
the \emph{support of $l$ with respect to $\gamma$ and $\Sigma,$} 
denoted by ${\rm Supp}(l),$ 
consists of the elements 
$v_j$ of $\Aa$ such that  $l_j + \gamma_j \notin \ZZ_{\geq 0}.$
We define the $\emph{set}$ 
$\Sa_{\Sigma}(\gamma) \subset \LL$
by the property that 
$l \in \LL$ belongs to 
$\Sa_{\Sigma}(\gamma)$ if there exists a (maximal) cone 
$\sigma$ such that all the elements of 
$\sup (l)$ generate rays of $\sigma.$
\end{definition}

Note that for any $l \in \LL$ and an arbitrary $\gamma \in \QQ^n,$
we have that 
$$
\Sa_\Sigma(\gamma)=(-l)+\Sa_\Sigma(\gamma-l).$$

\begin{remark}
For $v \in \Box(\Sigma)$ and $\gamma^v \in \QQ^n$ a corresponding 
solution to $\gamma^v_1 v_1 + \ldots + \gamma^v_n v_n= \beta,$
any $l \in \Sa_\Sigma (\gamma^v)$ has the property that 
the cone $\RR_{\geq 0} \sup (l)$  belongs to the fan $\Sigma,$ and 
$\sigma(v)$ is a subcone of the cone
$\RR_{\geq 0} \sup (l).$
\end{remark}

The motivation behind Definition \ref{def:coner} is the following
result.

\begin{proposition}\label{prop:zero}
The expression 
$$
x^{v} \prod_{j=1}^n  \frac{z_j^{l_j+ \gamma^v_j+D_j}}
{\Gamma(l_j+ \gamma^v_j+D_j+1)},
$$
vanishes in the completion of the ring $\CC[K, \Sigma],$
unless $l \in \Sa_{\Sigma} (\gamma^v).$
\end{proposition}

\begin{proof} 
Suppose that the expression is non-zero and let 
$\sigma(v)$ be the smallest cone of $\Sigma$ that contains $v.$

Notice that $j$-th factor in the product is divisible by $D_j$ 
for all negative integer $l_j+\gamma^v_j.$ 
This in particular implies 
that $l_j+ \gamma_j^v\geq 0$ for
$\RR_{\geq 0} v_j \not\in \Sigma$, since these $l_j+ \gamma^v_j$ are integer. 
Because of the factor
$x^v$ in the expression, the set of all $v_j$ such that
$l_j+\gamma^v_j$ is a negative integer 
must lie in a cone of $\Sigma$ that contains $\sigma(v).$ Then the
rays of this cone contain all $v_j$ for which $l_j+ \gamma^v_j$ 
is either negative or non-integer. 
\end{proof}

To any maximal cone $\sigma$ of the fan $\Sigma$ we 
associate the set 
$\Ca_{\sigma} \subset \LL \otimes \RR \subset \RR^n$ 
defined by the property that
$x=(x_1, \ldots, x_n) \in \LL \otimes \RR$ belongs to $\Ca_\sigma,$ if $x_j 
\geq 0$
whenever $\RR_{\geq 0} v_j$ is not a ray $\sigma.$ It follows that
$\Ca_\sigma$ 
is 
a cone in  $\RR^n$ generated by 
elements of $\LL.$ Moreover, the Minkowski sum
$$
\Ca_\Sigma:= \sum_{\sigma \in \Sigma(d)}  \Ca_{\sigma}
$$
is a cone in $\RR^n$ 
whose dual cone $\Ca_\Sigma^\dual \subset 
(\RR^n)^\dual$ has
non-empty interior (see, for example, page 219 in \cite{GKZbook}). 
The cone $\Ca_\Sigma^\dual$
is the maximal cone associated to the chosen regular triangulation 
in the non-pointed secondary fan determined in $(\RR^n)^\dual$
by the set $\Aa.$

\begin{lemma} \label{lemma:incl}
Let $\gamma=(\gamma_1, \ldots, \gamma_n)$ be an element
in $\QQ^n.$
There exists 
$b=(b_1, \ldots, b_n) \in \LL$ 
such that $\Sa_{\Sigma} (-b +\gamma) \subset \Ca_{\Sigma},$ 
or, equivalently, such that $\Sa_{\Sigma} (\gamma) 
\subset (-b)+ \Ca_{\Sigma}.$ 
\end{lemma}

\begin{proof}
Given $\gamma \in \QQ^n,$ 
for each maximal cone $\sigma \in \Sigma,$ there exists a unique
element $\gamma^{\sigma}=(\gamma^{\sigma}_1, \ldots, \gamma^{\sigma}_n) \in 
\LL \otimes \QQ$
such that $\gamma^{\sigma}_j = \gamma_j,$ when $\RR_{\geq 0}v_j \notin \sigma.$ 
Note that
the set of all $l=(l_1, \ldots, l_n) \in \LL$ 
such that $l_j +\gamma_j \geq 0$ when 
$\RR_{\geq 0}v_j \notin \sigma,$
is contained in the translated
cone $(-\gamma^{\sigma})+ \Ca_\sigma \subset \LL \otimes \RR,$ where 
$\Ca_\sigma$ is the cone described in the first part of the previous lemma. 
We have that $(-\gamma^{\sigma})+ \Ca_{\sigma} 
\subset (-\gamma^{\sigma})+ \Ca_{\Sigma}.$

We claim that there exists an element $b \in \LL$ such that
$(-\gamma^{\sigma}) + \Ca_{\Sigma} \subset (-b) + \Ca_{\Sigma},$ 
for all the maximal 
simplices $\sigma \in \Sigma.$ Indeed, it is enough to find 
$b \in \LL$ such that $b-\gamma^{\sigma} \in \Ca_{\Sigma}$ 
for all the
maximal simplices $\sigma \in \Sigma.$ This can be achieved by choosing
a finite number of generators in $\LL$ for the cone $\Ca_{\Sigma},$ and
choosing  
$b \in \LL$  to be 
an appropriate positive integral linear combination of the generators. 

It follows that
$$
\bigcup_{\sigma \in \Sigma(d)} (-\gamma^{\sigma})+ \Ca_\sigma
\subset (-b) + \Ca_{\Sigma},$$
where the union is taken over all the maximal cones $\sigma \in \Sigma.$
Hence
$$
\Sa_{\Sigma} (-b +\gamma) \subset
\bigcup_{\sigma \in \Sigma(d)} 
\big( (b-\gamma^{\sigma})+ \Ca_\sigma \big) \subset \Ca_{\Sigma}.
$$
This also means that
$$
\Sa_{\Sigma}(\gamma) =(-b)+ \Sa_{\Sigma} (-b + \gamma)
\subset (-b) + \Ca_\Sigma.$$
\end{proof}

\begin{corollary}\label{cor:choice}
For any $v \in \Box(\Sigma),$ with $v = \sum_{j=1}^n q^v_j v_j,$  
with $0 \leq q_j <1,$ and $q_j=0$ if $\RR_{\geq 0} 
v_j \not\prec \sigma(v),$
there exists a choice of  
$\gamma^v \in \QQ^n$ such that 
$$\gamma_1^v v_1 + \ldots + \gamma_n^v v_n = \beta, \ 
\gamma_j^v \equiv q_j^v \mod \ZZ,$$
with the property that
$\Sa_\Sigma (\gamma^v) \subset \Ca_\Sigma.$
\end{corollary}

\begin{proposition}\label{prop:cconv} 
For each maximal cone $\sigma$ of $\Sigma,$ 
there exists $c_\sigma \in {\mathcal C}_\sigma^\dual,$ such that the series 
$$
\sum_{l\in \LL \cap \Ca_{\sigma} }\;
\prod_{j=1}^n \frac{z^{\ l_j+ \lambda_j}_j} {\Gamma (l_j+\lambda_j+1)}.
$$
is absolutely convergent for $(z,\lambda) \in U_\sigma \times \CC^n,$ 
and defines an analytic function in $U_\sigma \times \CC^n.$
Here,
$
z_j^{\lambda_j}= e^{\lambda_j (\log |z_j| + i \arg z_j)},$
and the open set $U_\sigma$ 
is defined by
\begin{equation*}
\begin{split}
U_\sigma:= 
&\big\{ (z_1, \ldots, z_n) \in \CC^n : (-\log |z_1|, \ldots, -\log |z_n|) 
\in \Ca_{\sigma}^\vee + c_\sigma, \\ 
&(\arg z_1, \ldots, \arg z_n) \in (-\pi,\pi) \times
\ldots \times (-\pi,\pi) \big\}.
\end{split}
\end{equation*}
\end{proposition}

\begin{proof} The argument restriction on the $z$ variables
stems from the presence of the terms $z_j^{\lambda_j}.$
Hence, it is enough to show that the series
\begin{equation*}
\sum_{l\in \LL \cap \Ca_{\sigma} }\;
\prod_{j=1}^n \frac{z^{\ l_j}_j} {\Gamma (l_j+\lambda_j+1)}
\end{equation*}
converges absolutely for 
$(z,\lambda) \in \CC^n \times \CC^n \,$ with
$$
\big\{ z=(z_1, \ldots, z_n) \in \CC^n : (-\log |z_1|, \ldots, -\log |z_n|) 
\in \Ca_{\sigma}^\vee + c_\sigma
\big\},
$$
for some $c_\sigma \in \Ca_{\sigma}^\vee,$ to be determined below.

In order to be able to apply the Weierstrass convergence theorem
for sequences of holomorphic functions, we have to investigate 
the uniform convergence of the sequence of partial sums of this 
series for $\Vert \lambda \Vert \leq \delta.$ 

For any $l=(l_1, \ldots,l_n) \in \LL,$ we have $\sum_{j=1}^n l_j=0.$ 
Hence, we see that 
$$| \sum_{j=1}^n \Re (l_j+\lambda_j+1) | \leq \delta +n,$$ 
and 
$$\sum_{j=1}^n | \Im (l_j + \lambda_j +1)| \leq \delta +n,$$ 
so we can 
apply lemma \ref{lemma:est02} of the appendix. There exists a 
positive constant $A>0,$ such that
$$
\big| \prod_{j=1}^n \frac{z^{\ l_j}_j} {\Gamma (l_j+\lambda_j+1)} \big|
\leq A \cdot (4n)^{\Vert l \Vert} e^{\sum l_j \log |z_j|}=
A \cdot e^{\Vert l \Vert \log (4n) + \sum l_j \log |z_j|}.
$$

For some $\epsilon >0,$ the absolute and uniform convergence of the series 
is guaranteed for those $z=(z_1, \ldots, z_n) \in \CC^n$
such that
$$
\Vert l \Vert \log (4n) + \sum_{j=1}^n l_j \log |z_j| \leq -\epsilon 
\Vert l \Vert,
$$
for any $l \in \LL \cap \Ca_{\sigma}.$ 

Choose $\tilde{l}_1, \ldots, \tilde{l}_p$ to be a 
set of generators of the cone $\Ca_\sigma.$ We can 
choose $c_\sigma$ deep
enough in the interior of the 
cone $\Ca_{\sigma}^\dual,$ such that, for any $u \in 
\Ca_{\sigma}^\dual + c_{\sigma},$
$$
\langle u, \tilde{l}_i \rangle \geq (\epsilon + \log (4n) )
\log \Vert \tilde{l}_i \Vert $$
for any $i, 1 \leq i \leq p.$ This implies that,
for any $l \in \LL \cap \Ca_{\sigma}$ and 
any $u \in 
\Ca_{\sigma}^\dual + c_{\sigma},$ 
we have that 
$$
\langle u, l \rangle \geq (\epsilon + \log (4n) \rangle 
\log \Vert l \Vert.$$
It follows that the series 
converges absolutely in the region $U_\sigma$ introduced in 
the statement of the proposition. The 
region contains an open set of $\CC^n,$ since $\Ca_\Sigma^\vee
\subset \Ca_\sigma^\vee$ and 
the cone 
$\Ca_{\Sigma}^\vee$
has nonempty interior (page 219 in \cite{GKZbook}). This ends the proof
of the proposition. 
\end{proof}

\begin{corollary}\label{cor:conv1}
Let $J$ be a subset of the set of the maximal cones
of the fan $\Sigma,$ and 
$\Ca_J:= \sum_{\sigma \in J} \Ca_\sigma.$
Then there exists $c_J \in \Ca_J^\dual$ 
such that the series 
$$
\sum_{l\in \LL \cap \Ca_{J} }\;
\prod_{j=1}^n \frac{z^{\ l_j+ \lambda_j}_j} {\Gamma (l_j+\lambda_j+1)}.
$$
is absolutely convergent for $(z,\lambda) \in U_J \times \CC^n,$ 
and defines an analytic function in $U_J \times \CC^n,$
where
\begin{equation}\nonumber
\begin{split}
U_J:= 
&\big\{ (z_1, \ldots, z_n) \in \CC^n : (-\log |z_1|, \ldots, -\log |z_n|) 
\in \Ca_{J}^\vee + c_J, \\ &(\arg z_1, \ldots, \arg z_n) \in (-\pi,\pi) \times
\ldots \times (-\pi,\pi) \big\}.
\end{split}
\end{equation}
In the case $J=\Sigma(d),$ we will 
use the notations $c_\Sigma$ and $U_\Sigma$ to
denote $c_{\Sigma(d)}$ and $U_{\Sigma(d)},$ respectively.
\end{corollary}

\begin{proof} We have that $\Ca_\Sigma^\dual \subset 
\Ca_J^\dual=\cap_{\sigma \in J} \Ca_\sigma^\dual,$ therefore 
the cone $\Ca_J^\dual$ has nonempty interior. It is then 
enough to apply the previous proposition, and to note 
that, it is possible to choose $c_J \in \Ca_J^\vee$ such that 
$c_J -c_\sigma \in \Ca_J^\vee,$ for all $\sigma \in J.$
The corollary follows after we see that
$$
\Ca_J^\dual + c_J = \big(\Ca_J^\dual + (c_J - c_\sigma)\big)+ c_\sigma 
\subset \Ca_J^\dual + c_\sigma \subset \Ca_\sigma^\dual + c_\sigma,
$$
for all $\sigma \in J.$
\end{proof}

\begin{corollary}\label{cor:convs}
For any $v \in \Box(\Sigma),$ 
the series 
$$
\sum_{l\in \Sa_\Sigma(\gamma^v)}
\sum_{l\in \Sa_\Sigma(\gamma^v)}
\prod_{j=1}^n \frac{z_j^{l_j+ \gamma_j^v +
\lambda_j}} {\Gamma (l_j+\gamma_j^v+\lambda_j+1)}.
$$
defines an analytic function in the domain $(z,\lambda)
\in U_\Sigma \times \CC^n.$
\end{corollary}

\begin{proof}
By lemma \ref{lemma:incl}, we can find an element $b \in \LL$
such that $\Sa_\sigma (\gamma^v) \subset (-b) + \Ca_{\Sigma}.$
Hence, the above series is bounded in absolute value 
by the series
$$\sum_{l\in \LL \cap \Ca_\Sigma}\big\vert
\prod_{j=1}^n \frac{z_j^{l_j-b_j+\gamma_j^v +
\lambda_j}} {\Gamma (l_j-b_j+\gamma_j^v+\lambda_j+1)} \big\vert.
$$
The change of variable $\lambda_j \to \lambda_j -b_j+\gamma_j^v,$
implies that, by corollary \ref{cor:conv1} with $J=\Sigma(d),$ 
the series is absolute convergent for $(z,\lambda) \in U_\Sigma \times \CC^n.$
\end{proof}

\begin{definition}
Let $z=(z_1, \ldots, z_n)$ be a point in the open set
$U_\Sigma \subset \CC^n$ introduced in corollary \ref{cor:conv1}. 
The {\it $\Gamma$--series with values in the completion of  
$\CC[K,\Sigma]$}
is defined as 
$$
\Phi_{\Sigma}(z_1, \ldots, z_n):=\sum_{v \in \Box(\Sigma)} 
x^v \sum_{l\in \LL} 
\prod_{j=1}^n  \frac{z_j^{l_j+ \gamma^v_j+D_j}}
{\Gamma(l_j+ \gamma^v_j+D_j+1)},
$$
where, as before, 
$$
D_j := x^{v_j} \ \text{if} \ \RR_{\geq 0} v_j \in \Sigma, \ \text{and} \
D_j:=0, \ \text{otherwise,}$$ 
and 
$$
z_j^{D_j}:= e^{(\gamma^v_j+D_j) (\log |z_j| + i \arg z_j)},$$
for a choice of $(\arg z_1, \ldots, \arg z_n) \in \RR^n.$
\end{definition}

\begin{proposition}\label{prop:fct}
The series $\Phi_{\Sigma}(z_1, \ldots, z_n)$ defines 
a map from the region $U_{\Sigma} \subset \CC^n$ (defined in 
corollary \ref{cor:conv1})
to the completion of the graded ring $\CC[K,\Sigma].$
\end{proposition}

\begin{proof} According to proposition \ref{prop:zero},
for each $v \in \Box(\Sigma),$
the non-zero terms of the series come from $l \in
\Sa_\Sigma^{\gamma^v}.$ 
We then apply corollary \ref{cor:conv1} by setting $\lambda_j=
D_j,$ with $D_j = x^{v_j}$ if $\RR_{\geq 0} v_j \in \Sigma,$ and
$D_j=0,$ otherwise. The result follows.
\end{proof}

\subsection{GKZ solutions with values in SR--cohomology}

The equality 
$$
\Phi_{\Sigma}(z_1, \ldots, z_n)= \sum_{w \in K} \Phi_w (z_1,\ldots,z_n) x^w
$$
holds in the completion of the ring $\CC[K,\Sigma].$
Let $R$ be the subring of $\CC[K,\Sigma]$ 
generated by the elements $x^{v_j}$ for $\RR_{\geq 0} v_j\in \Sigma.$
We can view  $\CC[K,\Sigma]$ as an $R$-module.

\begin{definition}
The {\it leading term module} $M(\beta) \subset \CC[K, \Sigma]$ associated 
to the vector $\beta$ and  the fan $\Sigma$ 
is the $R$-submodule of $\CC[K, \Sigma]$
generated by the elements
$$
x^v \cdot \prod_{r_j<0, \RR_{\geq 0} v_j \in \Sigma, 
\RR_{\geq 0} v_j \not\prec \sigma(v)} x^{v_j},
$$
for all relations $v+ \sum r_j v_j=\beta,$
with $v \in \Box(\Sigma),$ $r \in \ZZ^n,$ such that
$r_j \geq 0$ if $\RR_{\geq 0} v_j\not\in\Sigma,$
and 
$\sigma(v)$ is the smallest cone that contains $v.$ 
\end{definition}

For a better understanding of the modules $M(\beta),$
let us choose $P > 0$ to be the least common multiple of all the indexes
of  sublattices in $\ZZ^d$ generated by all the possible
simplices with vertices among the vectors of $\Aa.$ In particular,
for any simplicial fan supported on the cone $K$ whose rays are 
generated by elements of $\Aa,$ $P v$ will be in the semigroup
generated by the elements of $\Aa,$ for any
$v$ in the twisted sector of that fan. We fix a 
fan $\Sigma.$

\begin{proposition} \label{prop:mb}
A lattice element $w \in N \cap K$ has the property that $x^w \in M(\beta)$ 
if and only if there exists some integer $k > 0$ such that
$$
(\beta - w) +  k P w
$$
is in the semigroup generated by all the elements of $\Aa.$
\end{proposition}

\begin{proof}
An element $w$ has the property that $x^w \in M(\beta)$ 
if and only if 
$$
w = v + \sum_{r_j<0, \RR_{\geq 0} v_j \not\prec \sigma(v)} v_j + \sum n_j v_j,
$$
for some $v$ in the twisted sector of the fan,
$n_j \in \ZZ_{\geq 0}$ and we only use $v_j$ generating rays of 
some maximal cone $\sigma$
of the fan that contains as a subcone the cone $\sigma(v),$
$$
v=\sum_{j=1}^n q_j v_j, \ 0 \leq q_j <1, q_j=0 \ 
\text{if} \ v_j \ \text{is not a ray of} \ \sigma(v).
$$
Here $r$ is a solution to $v + \sum_j r_j v_j =\beta,$ $r_j \in \ZZ,$
which
corresponds to the twisted sector $v,$
where the vectors $v_j$ corresponding to the 
negative  $r_j$ are rays of the maximal cone $\sigma.$
In particular, $r_j$ are nonnegative for $\RR_{\geq 0}v_j$ 
is not a cone of the fan.

Choose the positive integer $k$ such that
$r_j + kP q_j > 0,$ for $\RR_{\geq 0} v_j \prec \sigma(v),$ and
$r_j+kP -1>0$ for $\RR_{\geq 0} v_j \not\prec \sigma(v).$
Then
\begin{equation}\nonumber
\begin{split}
(\beta-w) + & kPw = \\
=&\sum_{\RR_{\geq 0}v_j \prec \sigma(v)} (r_j +kP q_j) v_j
+ \sum_{r_j<0, \RR_{\geq 0} v_j \not\prec \sigma(v)} (r_j+kP-1) v_j + \\
+&\sum_{r_j \geq 0} r_j v_j +\sum (kP-1)n_j v_j,
\end{split}
\end{equation}
which proves the only if part of the lemma.

To show the if part assume that for some integer $k>0$ we have
$$
(\beta-w) + kPw = \sum l_j v_j
$$
where $l_j \in \ZZ_{\geq 0}$ and the sum is taken over all the vectors $v_j$ 
in $\Aa.$ There exists a
maximal cone $\sigma$ and a corresponding twisted sector $v$
with the associated minimal cone $\sigma(v)$
such that
$$
w= v + \sum_{\RR_{\geq 0} v_j \prec \sigma} c_j v_j, c_j \in \ZZ_{\geq 0}.
$$
This allows us to write
$$
\beta= (1- kP) v - \sum_{v_j \in \sigma} (kP-1) c_j v_j
+\sum l_j v_j.
$$
We can see that in this presentation of $\beta,$ for those $v_j$
with $\RR_{\geq 0} v_j \not\prec
\sigma(v)$ that have negative coefficients, i.e. 
$-(kP-1)c_j + l_j <0,$ we must have $c_j >0.$
Then we can see that $w$ is in fact written as required in the
beginning of the proof.
\end{proof}

Note that the condition that $(\beta -w) + kPw$
belongs to the semigroup generated be the elements of $\Aa$ does not
depend on the choice of the simplicial fan $\Sigma.$ 

Recall that
$\CC[K^\circ, \Sigma]$ is the ideal of the ring $\CC[K, \Sigma]$
generated by the elements $x^w$ for $w \in K^\circ.$

\begin{corollary}\label{cor:bbb}
\

i) For any $\beta \in N,$ we have that 
$\CC[K^\circ, \Sigma] \subset M(\beta) \subset \CC[K, \Sigma].$

\

ii) If $\beta \in - K^\circ,$ then $M(\beta)= \CC[K^\circ, \Sigma].$

\end{corollary}

\begin{proof}
i) The second inclusion is true by the definition of the module
$M(\beta).$ 
To show the first inclusion observe that, for any  $w \in K^\circ,$
there exists a small rational number $\epsilon>0$
such that $w-\epsilon \sum v_j\in K$. Hence, $w$ can be written as a
\emph{positive} linear combination of all the elements 
$v_j$ with rational coefficients.
We use this representation of $w$ and any representation of $\beta-w$
as an integer linear combination of $v_j$ to see that for 
a sufficiently large and sufficiently divisible $k,$ the element
$\beta-w +kPw$ is a positive linear combination of $v_j$. Proposition
\ref{prop:mb} then completes the argument.

ii) It is enough to note that, when $\beta \in -K^\circ,$ the condition that $- \beta + (kP -1)w$ 
belongs to the semigroup generated by the elements of $\Aa$ implies that
$w \in K^\circ.$ 
\end{proof}

Recall that, for any basis $(m_1, \ldots, m_d)$ of ${\rm Hom}(N,\ZZ),$ 
the elements 
$$
Z_i=\sum_{j,\RR_{\geq 0} v_j \in \Sigma} \langle m_i, v_j \rangle
x^{v_j}
$$
form a regular sequence in $\CC[K,\Sigma].$ There is only a finite 
number of possible elements of
the form $\prod_{r_j < 0} x^{v_j},$ so the $R$-module
$M(\beta)$ is finitely 
generated. As a direct summand of the finite dimensional
vector space $\CC[K,\Sigma]/(Z_1, \ldots, Z_d)\CC[K,\Sigma],$
the vector space $R/(Z_1, \ldots, Z_d)\CC[K,\Sigma]$ is also
finite dimensional. Hence, the following holds:

\begin{proposition}
The vector space 
$$
M(\beta)/ZM(\beta):=M(\beta) / (Z_1, \ldots, Z_d) M(\beta)$$ 
is finite 
dimensional.
\end{proposition}

As a consequence of corollary \ref{cor:bbb}
and proposition \ref{prop:fct},
the $\CC[K, \Sigma]$--valued $\Gamma$-series $\Phi_{\Sigma}$
induces a map from
the region in $U_{\Sigma} \subset \CC^n$ to the finite dimensional vector
space $M(\beta) / Z M(\beta),$ 
$$
(z_1,...,z_n) \mapsto \Psi_{\Sigma} (z_1, \ldots, z_n):=
\sum_{v \in \Box(\Sigma)} \sum_{l\in \LL} 
R^{v}_{\gamma^v +l} (z_1,\ldots,z_n),
$$
where each term of the $\Gamma$--series   
is interpreted mod $(Z_1,\ldots,Z_d)M(\beta),$
or modulo $(Z_1,\ldots,Z_d) \CC[K, \Sigma],$ respectively.
By a slight abuse of notation, we will 
also denote by $\Psi_{\Sigma} (z_1, \ldots, z_n)$
the same map with values in the 
finite dimensional vector space $\CC[K,\Sigma] / Z \CC[K,\Sigma].$

\begin{proposition}
For any linear map $h : M(\beta)/ZM(\beta) \to \CC,$ 
(or,  $h : \CC[K,\Sigma] / Z \CC[K,\Sigma] \to \CC$) 
the function 
$h \cdot \Psi_{\Sigma}(z_1,\ldots,z_n)$ satisfies the GKZ hypergeometric equations
corresponding to the set $\mathcal A.$
\end{proposition}

\begin{proof}
The binomial GKZ equations
$$
\big(\prod_{l_j >0} \big(\frac{\partial}{\partial z_j}\big)^{l_j} -
\prod_{l_j <0} \big(\frac{\partial}{\partial z_j}\big)^{-l_j} \big)
\Psi_{\Sigma}=0,
$$
are satisfied because of the gamma identity.

For the linear GKZ equations, note first that, our choice that 
$D_j=x^{v_j}=0$ if $\RR_{\geq 0} v_j \notin \Sigma,$ shows that
$$
\sum_{j=1}^n \langle m, v_j \rangle D_j=
\sum_{v_j \in \Sigma} \langle m, v_j \rangle D_j,
$$
for any $m \in M={\rm Hom} (N, \ZZ).$ Hence
$$
\big( -\beta + \sum_{j=1}^n v_j z_j \frac{\partial}{\partial z_j} \big) \Psi_{\Sigma}
=\big(\sum_{j=1}^n v_j D_j \big) \Psi_{\Sigma}.
$$
It remains to observe that $\sum_{j=1}^n v_j D_j$ is a linear 
combination of the $Z_j$'s and that, by definition, the 
$\Psi_{\Sigma}(z_1, \ldots, z_n)$ takes values in $M(\beta) \subset 
\CC[K, \Sigma].$
\end{proof}

\begin{definition}\label{def:gammasr}
We call the induced maps 
$$
\Psi_{\Sigma} : U_{\Sigma} \to 
M(\beta) / Z M(\beta)$$
and 
$$\Psi_{\Sigma} : U_{\Sigma} \to 
\CC[K,\Sigma] / Z \CC[K,\Sigma]$$ 
the {\it GKZ solution map} and 
the {\it SR-cohomology valued GKZ solution map,} respectively.
\end{definition}

The next result deals with the linear independence of the solutions
obtained above.

\begin{proposition}
If $h : M(\beta)/ZM(\beta) \to \CC$ is a linear map
such that $h \cdot \Psi_{\Sigma}=0,$
then $h=0.$
\end{proposition}

\begin{proof} It is clear that solutions induced by elements of
$M(\beta)$ corresponding to different elements in
$\Box(\Sigma)$ are linearly independent, since 
the different fractional powers induce different monodromy 
behaviors.

Hence, we can restrict our attention to one element of  
$\Box(\Sigma)$ at a time. Consider $v \in \Box(\Sigma)$ with 
$$
v=\sum_{j=1}^n q_j v_j, \ 0 \leq q_j <1, q_j=0 \ \text{if} \ 
\RR_{\geq 0} v_j \not\prec \sigma(v).
$$

Assume that there exists an 
element $x \in M(\beta)/ZM(\beta)$ corresponding to
$v \in \Box(\Sigma)$ such that $h(x) \not= 0.$
Let $L$ be the largest degree of such an 
element. Furthermore, we can assume that 
$x$ lifts to a monomial, i.e. $x=x^w \mod ZM(\beta).$
Here
$$w= v+ \sum_{r_j<0, \RR_{\geq 0} 
v_j \not\prec \sigma(v)} v_j + \sum n_j v_j,$$ 
where
$\sum r_j v_j= \beta,$ $r_j \equiv q_j (\, {\rm mod} \, \ZZ),$ 
$n_j \in \ZZ_{\geq 0},$
and we only use $v_j$ generating rays of 
some maximal cone $\sigma$
of the fan that contains the cone $\sigma(v)$ as a subcone.
The 
non--negative integers $n_j$ are zero unless $v_j \in \Sigma.$

Let $\epsilon >0$ be some small positive number.
For each $j$ such that $v_j \in \Sigma,$ 
consider the loop of the form $z_j(t)= \epsilon e^{2 \pi i t}, 
z_i(t)=\epsilon,$
$i\not=j,$ $0 \leq t \leq 1.$ 
The action of the induced monodromy operator $T_j$
on the $\Gamma$-series $\Psi_{\Sigma}$ with values in $M(\beta)/ZM(\beta)$
is given by $exp(D_j)$. As a result, there is a polynomial
$g(T_j)$ such that $g(T_j) \Psi_{\Sigma} = D_j \Psi_{\Sigma}$, for every $j$
such that $\RR_{\geq 0} v_j \in \Sigma.$ As a result, we
have
$$
\prod_j g(T_j)^{n_j} h(\Psi_{\Sigma})(z)=
h (\prod_j D_j^{n_j} \Psi_{\Sigma}(z)).
$$
We now claim that the resulting function is nonzero.

The definition of the $\Gamma$--series  
$\Psi_{\Sigma}(z_1, \ldots, z_n)$ and the fact that $D_j=x^{v_j}$ are nilpotent in 
$M(\beta)/ZM(\beta)$ shows that any induced solution of the
GKZ system
can be written as the product of a monomial in the variables
$z_j$ and an element of $\CC[u_k^{-1}, \log u_k] [[u_k]]$
where $u_k, 1 \leq k \leq n-d,$ are torus invariant variables.
As a consequence of this fact, it is enough to show that the
``Fourier coefficient'' of $z^{r},$ for some $r,$ is nonzero.
We choose the element $r=(r_1, \ldots, r_n)$ introduced above,
with the property that $\sum r_j v_j= \beta.$

The ``Fourier coefficient'' of $z^r$ in the expansion of 
$h (\prod_j D_j^{n_j} \Psi_{\Sigma}(z))$ is given by
$$
h(\prod  D_j^{n_j} \cdot x^v \cdot \prod \frac{1}{\Gamma(r_j+D_j+1)}).
$$
Notice that the terms that occur in the expansion of the 
expression in the argument 
of $h$ have degree at least $L.$
Moreover, $x^w$ is the only monomial of that degree
that occurs, and it has a nonzero coefficient.
Since $h(x^w) \not= 0,$
the maximal property of $L$ implies that $z^r$ has  
a non-zero coefficient
times $h(x^w),$ so it is nonzero. This ends the proof of the 
linear independence result.
\end{proof}
The following result is proved by different methods in
section 3.5 of \cite{SST}.

\begin{proposition} For any $\beta \in N,$ 
$$
\dim M(\beta)/ZM(\beta) \geq {\rm Vol}(\Delta),
$$
where ${\rm Vol}(\Delta)$ is the normalized volume of the 
polytope $\Delta.$

\end{proposition}
\begin{proof}
We view $M(\beta)$ as a graded $\CC[y_1,...,y_d]$ module with $y_i=Z_i$.
Since $R$ is finitely generated as a module over this ring,
and since $M(\beta)$ is finitely generated over $R$, $M(\beta)$ is
finitely generated. By the graded Nakayama lemma, a basis of 
$M(\beta)/ZM(\beta)$ can be lifted to a set of generators of $M(\beta)$
as a module over $\CC[y_1,...,y_d]$.
By corollary \ref{cor:bbb}(i), the dimension of the $k$-th graded component
of $M(\beta)$ grows as a polynomial of the form 
$$
{\rm Vol}(\Delta) \frac {k^{d}}{(d-1)!} + {\rm lower~degree~terms.}
$$
Since the dimension of the degree $k$ component of $\CC[y_1,...,y_d]$ 
grows like $\frac {k^{d}}{(d-1)!}+...,$ the number of generators
of $M(\beta)$ is at least ${\rm Vol}(\Delta)$.
\end{proof}

\begin{corollary}\label{cor:gkzsvol}

\

i) For any $\beta \in N,$ the map
\begin{equation*}
\begin{split}
(M(\beta)/ZM(\beta))^\vee &\to {\mathcal Sol} (U_{\Sigma}) \\
f &\to f \cdot \Psi_{\Sigma}
\end{split}
\end{equation*}
produces at least
${\rm Vol} (\Delta)$ linearly independent GKZ solutions
which are analytic in $U_{\Sigma}.$

ii)
If $\, \beta \in -K^\circ,$ the above map
produces exactly ${\rm Vol} (\Delta)$ linearly independent GKZ solutions
which are analytic in $U_{\Sigma}.$
\end{corollary}

\section{GKZ solutions with values in $K$--theory}
\label{sec:gkzk}

We first recall some results about reduced toric 
Deligne-Mumford stacks and their $K$--theory. 
According to \cite{BCS},
a stacky fan $\Sigma$ in the abelian group $N$ is defined by a usual 
simplicial fan in $N \otimes \RR$ 
and a finite set of vectors $v_j$ ($1\leq j \leq n$) in $N$
generating the rays of the fan. In the context of this work, 
$N$ is always a lattice. 

According to \cite{BH},
the Grothendieck group $K_0 (\PP_\Sigma)$ is generated by the classes
$R_j$ of the invertible sheaves ${\mathcal L}_j$ corresponding to the 
one dimensional cones of the fan $\Sigma.$ Moreover, the ring $K_0(\PP_\Sigma)$
is isomorphic to the quotient of the Laurent polynomial
ring $\ZZ[R_1^\pm, \ldots, R_n^\pm]$
by the ideal generated by the relations:
\begin{itemize}
\item
$\prod_{j=1}^n R_j^{\langle m, v_j \rangle} 
=1,$ for all $m \in M= {\rm Hom} (N,\ZZ),$
\item
$\prod_{j \in J} (1-R_j)=0,$  for any set $J \subset 
\{1, \ldots, n \},$ such that
$\sum_{j \in J} \RR_{\geq 0} v_j$ is not a cone of the fan $\Sigma.$
\end{itemize}
The next proposition is a restatement of the main results of 
section 5 in \cite{BH}. 

\begin{proposition}\label{thm:kth}

\

i) The ring $K_0(\PP_\Sigma,\CC)$ is Artinian. Its maximum ideals
are in one-to-one correspondence with
elements of $\Box({\Sigma})$ as follows.
An element $v = \sum_{j=1}^n q^v_j v_j,$  
with $0 \leq q_j <1,$ and $q_j=0$ if 
$\RR_{\geq 0} v_j \not\prec \sigma(v),$
corresponds to the maximum ideal determined by 
the $n$-tuple of roots of unity
$(y^v_1,\ldots,y^v_n)\in\CC^n$ with $y^v_j=\ee^{2\pi i q^v_j}.$

ii) The  $K_0(\PP_\Sigma,\CC)$ is a direct sum 
of Artinian local rings obtained by localizing
at maximal ideal corresponding to all elements
$v \in \Box(\Sigma),$
$$
K_0(\PP_\Sigma,\CC) = \oplus_{v \in \Box(\Sigma)}
(K_0(\PP_\Sigma,\CC))_v.$$

iii) There is a natural vector space isomorphism between 
$K_0(\PP_\Sigma,\CC)$ and the SR--cohomology ring 
$\CC[K,\Sigma]/Z\CC[K, \Sigma].$ It is induced by isomorphisms
of $\CC$--algebras
\begin{equation*}
\begin{split}
(K_0(\PP_\Sigma,\CC))_v &\cong x^v \cdot \CC[K,\Sigma]/Z\CC[K, \Sigma] \\
&\cong H_{SR} (\PP_{\Sigma/\sigma(v)},\CC),
\end{split}
\end{equation*}
for each element $v \in \Box(\Sigma).$ Here  
$H_{SR} (\PP_{\Sigma/\sigma(v)},\CC)$ denotes the SR-cohomology ring of the 
toric smooth Deligne-Mumford stack
induced by the quotient fan $\Sigma/\sigma(v).$

iv) For any element $v \in \Box(\Sigma),$ with $v=\sum_{j=1}^n q^v_j v_j$
with $0 \le q^v_j < 1,$ the $\CC$-algebra isomorphism 
$$
(K_0(\PP_\Sigma,\CC))_v \cong 
H_{SR} (\PP_{\Sigma/\sigma(v)},\CC)$$
is given by 
$$
R_j=y^v_j e^{D_j}, \ 1 \leq j \leq n,
$$
where $y^v_j=e^{2 \pi i q^v_j},$
and $R_j$ and $D_j$ are the generators of the two rings
corresponding to the vectors $v_j.$ 
\end{proposition}

A few facts and some notation used in the spectral theory of linear operators
on finite dimensional vector spaces are collected in the Appendix \ref{app:2}. 
We will use them to construct GKZ solutions associated to the set 
$\Aa.$ 

Assume that the stacky 
fan $\Sigma$ is supported on the cone $K$ generated
by the elements of the set $\Aa,$ and that $v_j \in \Aa.$
In particular, consistent with 
the second set relations that hold 
in $K_0(\PP_{\Sigma}),$ we can  
assume $R_j=1$ 
whenever $\RR_{\geq 0} v_j$ is not a cone in $\Sigma.$
The linear operators 
$\Ra_j : 
K_0(\PP_\Sigma,\CC) \to K_0(\PP_\Sigma,\CC),$ $1 \leq j \leq n,$
defined by $\Ra_j (x):= R_j x,$ are mutually commuting with spectra
$s(\Ra_j)= \{ y_j^v \ : \ v \in \Box(\Sigma) \}$ corresponding
to the direct sum decomposition of $K_0 (\PP_{\Sigma}, \CC).$
We cover each root of unity $y^v_j$ ($v \in \Box(\Sigma), 1\leq j \leq n$)
with an open disc $B_j^v$ of radius $\epsilon > 0,$ such that any two discs 
centered at different points have disjoint closures, and  
the origin of the complex plane is not inside any of them.
The multidiscs
$
B^v:= B_1^v \times \ldots \times B_n^v \subset \CC^n
$
have the property that $\bar{B}^v \cap \bar{B}^{v^\prime} = \emptyset,$
for any two disjoint $v,v^\prime$ in $\Box(\Sigma).$ 

For an element $v \in {\rm Box} (\Sigma),$ 
and $U \subset \CC^n$ an open simply connected domain, 
consider the analytic function $\varphi : U \times B^v \to \CC,$
$$
\varphi (z, r):= 
\prod_{j=1}^n 
\frac{z_j^{\frac{1}{2\pi i} \log_j r_j}}
{\Gamma(\frac{1}{2\pi i} \log_j r_j+1)},
$$
where $\log_j$ are arbitrary log branches.
The spectrum of the restriction of operator $\Ra_j$ to 
$(K_0(\PP_\Sigma,\CC))_v$ consists of 
the unique 
value $y^v_j \in B^v_j.$ Since the function $\varphi (z,r)$ is
analytic for $z\in U$ and $r=(r_1, \ldots, r_n) \in B^v_1 \times \ldots \times B^v_n,$
the linear operator 
$$\varphi (z, \Ra) : (K_0(\PP_\Sigma,\CC))_v \to (K_0(\PP_\Sigma,\CC))_v$$
is well defined.

\begin{proposition}\label{prop:kzero}
The linear operator 
$$\varphi (z, \Ra) : (K_0(\PP_\Sigma,\CC))_v \to (K_0(\PP_\Sigma,\CC))_v$$
is zero, unless there exists   
a cone $\sigma \in \Sigma$ 
with $\sigma(v) \subset \sigma$
whose rays 
contain all $v_j \in \Aa$ for which
$
\frac{1}{2 \pi i}\log_j y^v_j 
$
is not a nonnegative integer. 
\end{proposition}

\begin{proof} 
The result is just a rephrasing of 
proposition \ref{prop:zero}. 
For $r_j \in B_j^v,$ we can write that 
$$
\frac{1}{2\pi i} \log_j r_j = \frac{1}{2\pi i}
\log_j  y^v_j + (r_j - y^v_j) \psi(r_j),
$$
where $\psi$ is analytic in $B^v_j.$ It follows that, 
in the domain $B^v,$ the function
$\varphi(z,r)$ is the product of an analytic function and the analytic function
$
\prod_{j \in J} (r_j - 1),$
where $J \subset \{ 1, \ldots, n \}$ consists of all $j$ such that
$
\frac{1}{2 \pi i}\log_j  y_j^v 
$
is a negative integer. 

Theorem VII.1.5(b) in \cite{DS} implies that the linear operator 
$\varphi(z,\Ra)$ is then
the product of a linear operator and the linear operator
$$
\prod_{j \in J} (\Ra_j -I).$$
According to
theorem \ref{thm:kth} $iv),$ the action of this operator on 
the space $(K_0(\PP_\Sigma,\CC))_v$ is a multiple of the 
element $\prod_{j \in J} D_j$ viewed as an element in the
SR cohomology ring $H_{SR}(\PP_{\Sigma/\sigma(v)}, \CC).$ The
combinatorial definition of this ring provides the final step 
of the proof.
\end{proof}

For some $v \in \Box(\Sigma),$ we make 
a choice of log branches, such that, 
$\gamma_1^v v_1 + \ldots + \gamma_n^v v_n = \beta,$ 
where
$
\gamma_j^v:= \frac{1}{2 \pi i} \log_j y^v_j.$
Corollary \ref{cor:convs} shows that the function 
$\Xi_\Sigma^v (z,r):  U_{\Sigma} \times B^v \to \CC,$
defined by
\begin{equation}\label{eq:psiv}
\Xi_\Sigma^v (z,r):=
\sum_{l \in \Sa_{\Sigma}(\gamma^v)} 
\prod_{j=1}^n \frac{z_j^{l_j+\frac{1}{2\pi i} \log_j r_j}}
{\Gamma(l_j+\frac{1}{2\pi i} \log_j r_j+1)}
\end{equation}
is analytic.

\begin{corollary}
For some $v \in \Box(\Sigma),$ and $z \in U_\Sigma,$ 
the linear operator $\Xi_\Sigma^v (z,\Ra) : 
(K_0(\PP_\Sigma,\CC))_v \to (K_0(\PP_\Sigma,\CC))_v$
is given by 
$$
\Xi_\Sigma^v (z,\Ra)=
\sum
\prod_{j=1}^n \frac{z_j^{\frac{1}{2\pi i} \log_j \Ra_j}}
{\Gamma(\frac{1}{2\pi i} \log_j \Ra_j+1)},
$$
where the summation is taken over all the possible
$\log$ branches such that
$$
(\log_1 y^v_1)v_1 + \ldots + (\log_n y^v_n)v_n= 
(2 \pi i) \beta.$$
\end{corollary}

\begin{proof} 
For a fixed choice of log branches as above, and 
$
\gamma_j^v= \frac{1}{2 \pi i} \log_j y^v_j,$
proposition \ref{prop:kzero} implies that the nonzero terms
in the summation
$$
\sum_{l \in \LL}
\prod_{j=1}^n \frac{z_j^{l_j+ \frac{1}{2\pi i} \log_j \Ra_j}}
{\Gamma(l_j+\frac{1}{2\pi i} \log_j \Ra_j+1)},
$$
correspond to those $l \in \LL$ with $l \in \Sa_{\Sigma} (\gamma^v).$ 
\end{proof}

\begin{definition}\label{def:kgamma}
The {\it $\Gamma$--series with values in $K$--theory} 
is the map $\Xi_{\Sigma}(z, \Ra)$
from the region $U_\Sigma \subset \CC^n$ to 
the space of linear operators on 
$K_0(\PP_\Sigma,\CC).$ For any 
$z \in U_\Sigma,$ 
the linear operator $\Xi_{\Sigma}(z, \Ra): 
K_0(\PP_\Sigma,\CC) \to  K_0(\PP_\Sigma,\CC)$ is associated
to the analytic function
$$
\Xi_{\Sigma}(z,r) : U_{\Sigma} \times B
\to \CC
$$
with the property that
$
{\Xi_{\Sigma}} (z,r)= \Xi_{\Sigma}^v (z,r),$ for all $v \in \Box(\Sigma),$ 
and $(z,r) \in  U_{\Sigma} \times B^v,$
where $B$ is a domain in $\CC^n$ such that
$$
\bigcup_{v \in \Box(\Sigma)} B^v \subset B.
$$ 
\end{definition}

\begin{remark}\label{rem:gene}
The definition allows for some ambiguity in 
in choosing the domain $B$
and the function $\Xi_{\Sigma}(z, r).$
However, the operator $\Xi_{\Sigma}(z, \Ra)$ is 
independent of these choices.
\end{remark}

\begin{remark} 
If we regard the $K$--theory ring $A= K_0(\PP_\Sigma,\CC)$
as a module over itself, we see that
the linear operators $\Ra_j :  A \to A,$
given by $\Ra_j(x)=R_j x$ are $A$--module endomorphisms.
For an arbitrary commutative ring $A,$ the map $e: {\rm End}_A (A) \to A,$ defined by
$
e(\Phi):=\Phi(1)$
is a ring isomorphism, so we can regard 
the $K$--theory endomorphisms as $K$--theory elements. 
\end{remark}

\begin{proposition}
For any $\beta \in N,$ 
the $\Gamma$-series with values in SR-cohomology (cf. definition 
\ref{def:gammasr})
$$\Psi_{\Sigma} (z) : U_{\Sigma}
\to \CC[K, \Sigma]/Z\CC[K,\Sigma]$$ can be written as
$$
\Psi_{\Sigma}(z)=Ch \big( \Xi_{\Sigma} (z,\Ra)(1)\big),$$ 
where 
$Ch : K_0 (\PP_\Sigma,\CC) \to \CC[K, \Sigma]/Z\CC[K,\Sigma]$
is the $\CC$-algebra isomorphism between the $K$--theory and 
the SR-cohomology ring of proposition \ref{thm:kth} \, iii),
and  $\Xi_{\Sigma} (z,\Ra)(1): U_{\Sigma} \to  K_0 (\PP_\Sigma,\CC)$ is the 
is $K$--theory valued $\Gamma$-series.
\end{proposition}

\begin{proof} It is clear that it is 
enough to check the statement for each twisted 
sector $v \in {\rm \Box}(\Sigma).$ 
We have to prove that 
$$
\Psi_{\Sigma}(z)=Ch \big( \Xi_{\Sigma}^v (z,\Ra)(1)\big)$$ 
in $(K_0(\PP_\Sigma, \CC))_v \cong H_{SR} (\PP_{\Sigma/\sigma(v)},\CC).$

We have that
$\Xi_{\Sigma}^v (z,\Ra)(1)$ is equal in 
$(K_0(\PP_\Sigma, \CC))_v$ to
$$
\sum_{0\leq j_1, \ldots j_n< \nu}
\partial_1^{j_1} \ldots \partial_n^{j_n}
\Xi_{\Sigma}^v(y_1^v, \ldots, y_n^v) \, 
\frac{(R_1-y_1^v)^{j_1}}{j_1 !} \ldots 
 \frac{(R_n-y_n^v)^{j_n}}{j_n !}
$$ 
for a sufficiently large positive integer $\nu.$ It follows that, 
as an element of $H_{SR} (\PP_{\Sigma/\sigma(v)},\CC),$
$Ch \big( \Xi_\Sigma^v (z,\Ra)(1)\big)$ 
is equal to 
$$
\sum_{0\leq j_1, \ldots j_n< \nu}
\partial_1^{j_1} \ldots \partial_n^{j_n}
\Xi_{\Sigma}^v(y_1^v, \ldots, y_n^v) \, 
\frac{(y^v_1 e^{D_1}-y_1^v)^{j_1}}{j_1 !} \ldots 
 \frac{(y^v_n e^{D_n}-y_n^v)^{j_n}}{j_n !}.
$$
This is the Taylor polynomial approximation of the 
analytic function 
$\Xi_{\Sigma}^v (z, (y_1^v e^{D_1}, \ldots, y_n^v e^{D_n}))$
around the origin $D_1=\ldots=D_n=0.$ Note that
$$
\Xi_{\Sigma}^v (z, (y_1^v e^{D_1}, \ldots, y_n^v e^{D_n}))=
\sum_{l\in \Sa_{\Sigma}(\gamma^v)} 
\prod_{j=1}^n \frac{z^{l_j+ \gamma_j^v +D_j}} {\Gamma (l_j+ 
\gamma_j^v +D_j + 1)},
$$
as analytic functions with $z \in U_{\Sigma}$ and $(D_1, \ldots, D_n)$ in a 
neighborhood of the origin in $\CC^n.$ Hence, 
the above Taylor polynomial is equal to $\Psi_{\Sigma}(z)$
when viewed as elements of $H_{SR} (\PP_{\Sigma/\sigma(v)},\CC).$
\end{proof}

\begin{corollary}\label{cor:msmap}
For any $\beta \in N,$ the mirror symmetry map
$MS_{\Sigma}:(K_0 (\PP_\Sigma,\CC))^\vee  
\to {\mathcal Sol} (U_{\Sigma})$
given by 
$$
MS_{\Sigma}(f):= f \big(\Xi_{\Sigma} (z,\Ra) (1)\big)
$$
produces 
GKZ solutions which are analytic in $U_{\Sigma}.$
\end{corollary}

\begin{remark}
The dimension of the space of GKZ solutions with 
values in $K$--theory or $SR$-cohomology
is generally not easy to calculate or even estimate.
Moreover, the $K$--theoretic meaning of the leading
term module $M(\beta),$ for general $\beta,$ is unclear but 
it is perhaps worth investigating.
\end{remark}

\section{Analytic Continuation of Hypergeometric Series}
\label{sec:achs}
Consider an {\it oriented} 
edge of the secondary polytope starting at the vertex of the secondary 
polytope corresponding to some regular triangulation $\Ta_+$ 
and ending at the vertex corresponding to another regular triangulation
$\Ta_-,$ and let $\Sigma_+$ and $\Sigma_-$  be the induced fans
supported on the cone $K.$
{\it In order to ease up some of the heavy notation, in this 
section and the next one, we will 
usually replace the subscript $\Sigma_\pm$ by $\pm.$ For 
example, we will 
write $\Sa_{\pm} (\gamma)$ instead of $\Sa_{\Sigma_{\pm}} (\gamma).$}

Theorem 2.10 in \cite{GKZbook}
shows that there exists a circuit $I$ (i.e a minimal linearly
dependent set of elements) in $\mathcal A$ determining an integral
relation of the form 
$$
h_1 v_1 + \ldots + h_n v_n = 0, \ h=(h_1, \ldots, h_n)\in \LL, 
$$
with $I= I_+ \cup I_-,$  where 
$$
I_+:= \{ v_j : h_j > 0 \}, \ I_-: = \{v_j: h_j < 0 \},
$$
such that the triangulations are both supported on the circuit $I$ 
and are obtained by a modification based on the circuit $I$
(see pages 231-233 in the book \cite{GKZbook} 
for detailed explanations). Moreover,
each of the fans $\Sigma_\pm$ has the property that, for 
every subset $\Fa \subset \Aa \setminus I,$ if $\Fa \cup 
(I \setminus v_0) $ are the rays of 
a maximal cone in $\Sigma_\pm$
for {\it some} $v_0 \in I_\pm,$ then the elements of  
$ \Fa \cup (I \setminus v)$ are the rays of a maximal cone in 
$\Sigma_\pm$ for {\it any} $v \in I_\pm.$
Such a subset $\Fa$ is said to be a {\it separating} set for the 
fans $\Sigma_\pm.$ 
Furthermore, the modification is obtained by replacing the set of  
all maximal cones generated by sets of the form 
$\Fa \cup (I \setminus v)$
($v \in I_+$) of $\Sigma_+$ 
with the set of all maximal cones generated   
by $ \Fa \cup (I \setminus v)$ ($v \in I_-$),
with $\Fa$ a separating set. 

\begin{definition}\label{def:ess}
We say that a maximal cone of the fans $\Sigma_\pm$ generated 
by a set of the form $\Fa \cup (I \setminus v),$
with $\Fa$ a separating set and $v \in I_\pm,$ 
is {\it essential}. We denote by $\Sigma^{es}_\pm(d)$
the sets of essential maximal cones of the fans $\Sigma_\pm.$
\end{definition}

We assume that 
the element $h=(h_1, \ldots, h_n)\in \LL$ is {\it primitive} in $\LL,$
i.e. it is not a non-trivial integral multiple of any other element 
in the lattice. This means that $\LL / \langle h \rangle$
is itself a lattice. We introduce the notation
$$
\LL^\prime: =  \LL/ \langle h \rangle,
$$
and let $p: \LL \to \LL^\prime$ denote the canonical projection. 

\begin{definition}\label{def:ir}
For any $r=(r_1, \ldots, r_n) \in (\CC^*)^n,$ we denote by 
$\Ia(r) \subset \CC^*$ the finite set of complex numbers 
$t$ 
such that $r_j t^{h_j}=1$ for some $j$ with $v_j \in I_-.$ 
\end{definition}

\begin{remark}
For a given $r=(r_1,\ldots, r_n) \in (\CC^*)^n,$ note that 
any two values 
$t,s \in \Ia(r)$ such that 
$(r_1t^{h_1}, \ldots, r_n t^{h_n})=
(r_1s^{h_1}, \ldots, r_n s^{h_n}),$ are in fact equal. 
Indeed, since $r_j$ are nonzero, we have that 
$t^{h_j}=s^{h_j}$ for all $j.$ But $(h_1, \ldots, h_n) \in \LL$
is primitive, so 
we can find integers $d_1, \ldots, d_n$ such that
$d_1 h_1 + \ldots + d_n h_n=1.$ Then $t=t^{d_1 h_1 + \ldots + d_n h_n}=
s^{d_1 h_1 + \ldots + d_n h_n}=s.$
\end{remark}

Let $\Box (\Sigma^{es}_{\pm}) \subset \Box (\Sigma_{\pm})$
be the subsets
consisting of those elements $v \in \Box (\Sigma_{\pm})$ 
with the property that 
the minimal cones containing $v$ in $\Sigma_\pm$
are subcones of one of the maximal 
cones in $\Sigma_{\pm}^{es}(d),$ 
respectively. We now describe the effect that a 
modification has on the corresponding twisted sectors. We use
the notation and the results of Proposition \ref{thm:kth} $i).$

\begin{proposition}\label{prop:boxx}
\

i) $$\Box (\Sigma_{+}) \setminus \Box (\Sigma^{es}_{+})=
\Box (\Sigma_{-}) \setminus \Box (\Sigma^{es}_{-}).
$$
Moreover, for any element $v \in N$ belonging to the two sets
above, with
$$
v=\sum_{j=1}^n q_j v_j, \ 0 \leq q_j <1, q_j=0 \ \text{if} \ 
\RR_{\geq 0} v_j \not\prec \sigma(v),
$$
the minimal cone $\sigma(v)$ containing $v$ is unchanged
under the modification. 
Moreover, the corresponding sets of $n$-tuples of roots of unity 
$(y_1^v, \ldots, y_n^v),$ $y_j^v= e^{2 \pi i q_j},$ are 
also unchanged. 

ii) For any $v \in \Box (\Sigma^{es}_{+}),$ with 
$$
v=\sum_{j=1}^n q_j v_j, \ 0 \leq q_j <1, q_j=0 \ \text{if} \ 
\RR_{\geq 0} v_j \not\prec \sigma(v) \in \Sigma_+,
$$
and $y_j^v= e^{2 \pi i q_j^v},$ 
the  $n$-tuples $(y_1^v t^{h_1}, \ldots, y_n^v t^{h_n}) \in (\CC^*)^n$
with $t \in \Ia(y^v),$
determine maximum ideals of $K_0(\PP_{\Sigma_-}, \CC)$
corresponding to 
elements of $\Box (\Sigma^{es}_{-}).$ Moreover, any
$n$-tuple of roots of unity induced by some
element of $\Box (\Sigma^{es}_{-})$ is obtained
in this way, for appropriate $v \in \Box (\Sigma^{es}_{+})$
and $t \in \Ia(y^v).$
\end{proposition}

\begin{proof} $i)$ For any $v \in \Box (\Sigma_{+}) \setminus 
\Box (\Sigma^{es}_{+}),$
we have that $\sigma(v)$ is a subcone of a maximal cone of $\Sigma_+$
that survives the modification. Hence $\sigma(v)$ is also a cone
of $\Sigma_-.$ We still have to show that $v \notin \Box(\Sigma^{es}_-).$
For, assume that $\sigma(v)$ is a subcone of maximal 
cone $\sigma \in \Sigma_-^{es}(d)$ that changes under the modification,
and the rays of $\sigma$ consist of the elements of  
$\Fa \cup (I \setminus v_-)$
for some separating set $\Fa$ and $v_- \in I_-.$ Since
$\sigma(v)$ is a cone of $\Sigma_+$ and 
the elements of $I_+$ do not generate a cone in this fan,
there exists some vector $v_+ \in I_+$ which does not generate 
a ray of $\sigma(v).$ 
Hence $\Fa \cup (I \setminus v_+)$ is a 
maximal cone in $\Sigma^{es}_{+}(d)$ containing $\sigma(v)$
as a subcone. But this is a contradiction with our
starting assumption that $v \notin \Box (\Sigma^{es}_{+}).$
The roles of $\Sigma_+$ and $\Sigma_-$ can be reversed,
and the statement follows.

$ii)$ 
First, consider 
some $n$-tuple $(y_1^v, \ldots, y_n^v)$ corresponding to an 
element $v \in \Box(\Sigma_+^{es}).$ We show that 
the $n$-tuple 
$(y_1^v t^{h_1}, \ldots, y_n^v t^{h_n})$ with $t \in \CC$
a root of unity such that  
$y_i^v t^{h_i} =1,$ for some $i$ with $v_i \in I_-,$
corresponds to an element of $\Box(\Sigma_-^{es}).$
For, note that there exists a maximal cone $\sigma \in \Sigma_+^{es}(d)$
generated by $\Fa \cup (I \setminus w)$
for some separating set $\Fa$ and $w \in I_+,$
such that
$$v=\sum_{j=1}^n q_j v_j, \ 0 \leq q_j <1, q_j=0 \ \text{if} \ 
\RR_{\geq 0} v_j \not\prec \sigma.
$$
Choose a rational number $q \in \QQ$ such that
$q_{i} + q h_{i} \in \ZZ,$ for some $i$ with $v_i \in I_-.$
By adding the linear relation
$qh_1 v_1+ \ldots +qh_n v_n=0$ to the above expression of $v,$
we can write that
$$
v=\sum_{j=1}^n (q_j + q h_j) v_j, 
$$
where $q_j + q h_j \in \ZZ,$ whenever 
$v_j \notin \Fa \cup (I \setminus v_{i}).$
Since the elements of 
$\Fa \cup (I \setminus v_{i})$ generate the rays of a maximal cone
in $\Sigma_-^s(d),$ we conclude that 
$(y_1^v t^{h_1}, \ldots, y_n^v t^{h_n}),$ with
$y_j^v= e^{2 \pi i q_j^v}, t=e^{2 \pi i q},$
corresponds indeed to an element in $\Box(\Sigma_-^{es})$ which differs
from $v$ by an integral linear combination of the vectors $v_j \in \Aa.$
This shows that, for any $v \in \Box(\Sigma_+^{es}),$
the procedure described in the second part of the 
proposition produces a subset of $\Box(\Sigma_-^{es}).$ 
We still have
to show that the union of all these subsets, for all $v 
\in \Box(\Sigma_+^{es}),$ is equal to $\Box(\Sigma_-^{es}).$

It is enough to show that, after starting
with an $n$-tuple of roots of unity 
$(y_1^v, \ldots, y_n^v)$ corresponding to
some element $v \in \Box(\Sigma_+^{es})$ and applying 
the above procedure from $\Box(\Sigma_+^{es})$ to $\Box(\Sigma_-^{es})$ 
and back, the $n$-tuple $(y_1^v, \ldots, y_n^v)$
is recovered. Note that for such an $n$-tuple there 
exists an element $v_i \in I_+$ such that 
$y_i^v=1.$ It follows that, 
if $(y_1^v t^{h_1}, \ldots, y_n^v t^{h_n})$
(for the appropriate $t \in \CC$)
corresponds to an element in $\Box(\Sigma_-^{es}),$
then $(y_1^v t^{h_1} (t^{-1})^{h_1}, \ldots, 
y_n^v t^{h_n} (t^{-1})^{h_n})$ in $\Box(\Sigma_+^{es})$
is an allowed choice under the procedure, 
since
$y_i^v t^{h_i} (t^{-1})^{h_i}=y_i^v=1,$ with $v_i \in I_+,$
as noted above.
This ends the proof of the proposition.
\end{proof}

For an element $v \in \Box(\Sigma_+),$ consider 
the associated 
$y^v \in (\CC^*)^n.$ Choose
$n$ $\log$ branches, all of which are
denoted $\log_+$ by a slight abuse of notation, 
such that
$\gamma^v_1 v_1+ \ldots + \gamma^v_n v_n= \beta,$
where
$$\gamma^v_j:=\frac{1}{2 \pi i} \log_+ y_j^v.$$
We now choose a branch of $\log t,$ and for 
any $t \in \Ia(y^v),$ we set 
$$
\gamma^v_j(t):= \gamma_j^v + h_j \log t.$$
We still have that 
$\gamma^v_1(t) v_1+ \ldots + \gamma^v_n(t) v_n= \beta.$
For each $t \in \Ia(y^v),$
it is possible to find $\log$ branches that 
will be denoted $\log_-$
such that 
$$\gamma^v_j(t)=\frac{1}{2 \pi i} \log_- (y_j^v t^{h_j}).$$

We now exhibit some special subsets 
of the sets 
$\Sa_{\pm} (\gamma) \subset \LL$
introduced in definition \ref{def:coner}. 

\begin{definition}\label{def:coners}
For an arbitrary $\gamma=(\gamma_1, \ldots, \gamma_n) \in \QQ^n,$ 
we define the $\emph{set}$ 
$\Sa^{es}_{\pm}(\gamma^v) \subset \LL$
by the property that 
$l \in \LL$ belongs to 
$\Sa^{es}_{\pm} (\gamma)$ 
if there exists a maximal cone $\sigma \in \Sigma_{\pm}^{es}(d)$
such that all the elements of $\sup (l)$ with respect 
to $\Sigma_{\pm}$ and $\gamma$ generate rays of $\sigma.$
\end{definition}

When there is no danger of confusion, we will simply 
talk about $\sup (l)$ with no mention of the fan and 
the element $\gamma.$

\begin{proposition}\label{prop:sss}
For any $v \in \Box (\Sigma_+),$ with $\gamma^v, \gamma^v(t)$ 
($t \in \Ia(y^v)$) chosen as above, the following are true:

i) If  $1 \notin \Ia(y^v),$ then
$\Sa^{es}_{+} (\gamma^v)= \Sa_{+} (\gamma^v).$ 
If $1 \in  \Ia(y^v),$ 
then $\gamma^v(1)=\gamma^v,$ and 
$$
 \Sa_{+} (\gamma^v) \setminus 
\Sa^{es}_{+} (\gamma^v)=  \Sa_{-} (\gamma^v)
\setminus \Sa^{es}_{-} (\gamma^v).
$$

ii) Under the natural projection $p: \LL \to \LL^\prime,$
the images of the sets $\Sa^{es}_{+} (\gamma^v)$  and 
$\Sa^{es}_{-} (\gamma^v(t))$ coincide, for any $t \in \Ia(y^v).$
\end{proposition}

\begin{proof} $i)$ Assume that there exists $l \in 
\Sa_{+} (\gamma^v) \setminus 
\Sa^{es}_{+}(\gamma^v).$ This means
that there exists a maximal
simplex $\sigma \in \Sigma_+$ such that 
all the vectors in $\sup (l)$ generate rays in 
$\sigma.$ 
Since any such simplex $\sigma$ is unchanged under
the modification, and the rays of $\sigma(v)$ are in 
the support of $l,$
we see that $v \in \Box(\Sigma_-)$ 
since $\sigma \in \Sigma_-(d).$ 
Hence $l \in \Sa_{-}(\gamma^v)$ and $1 \in \Ia(y^v)$
with $\gamma^v(1)=\gamma^v.$

We still have to show that, if $l \in 
\Sa_{+} (\gamma^v) \setminus 
\Sa^{es}_{+} (\gamma^v),$ 
then $l \notin \Sa^{es}_{-} (\gamma^v).$
For, assume that there exists a maximal cone $\sigma^\prime \in 
\Sigma_-^{es}(d),$ $\sigma^\prime$ whose rays are
generated by the elements of 
$\Fa \cup (I \setminus v),$ $v \in I_-,$   
such that the elements of the support of $l$ 
generate rays in $\sigma.$
Note that,
since $l \in \Sa_{+} (\gamma^v)$ and the elements
of $I_+$ do not generate a cone in $\Sigma_+,$
there must be 
some $i$ with $v_i \in I_+$ such that $l_i + \gamma_i^v 
\in \ZZ_{\geq 0}.$ This means that the elements 
of the support of $l$ generate rays of the cone 
of $\Sigma_+^{es}(d)$ generated by 
$\Fa \cup (I \setminus v_i).$  However, this 
contradicts the assumption that $l \notin 
\Sa^{es}_{+}(\gamma^v).$ With this, we have 
proved that 
$\Sa_{+} (\gamma^v) \setminus 
\Sa^{es}_{_+} (\gamma^v) \subset  
\Sa_{-} (\gamma^v)
\setminus \Sa^{es}_{-} (\gamma^v).$ 
The inverse inclusion is obtained by employing the 
completely analogous argument with the roles of  
$\Sigma_\pm$ reversed. 

$ii)$ Given $l \in \LL,$ the elements of the 
fiber $p^{-1}(p(l))$ consist of elements 
of the form $l +mh$ with $m \in \ZZ.$

Let $l$ be an element of $\Sa^{es}_{+} (\gamma^v).$ 
This means that the vectors $v_j$ such that 
$l_j + \gamma^v_j \notin \ZZ_{\geq 0}$
are among the elements of a set of the form
$\Fa \cup (I \setminus v_i)$ with $\Fa$ a 
separating set and $v_i \in I_+.$ In particular,
we see that the minimal cone of 
$\Sigma_+$ that contains $v$ 
is itself contained in a maximal cone of $\Sigma_+$
that gets replaced under the 
modification. Hence $v \in \Box(\Sigma^{es}_+).$
According to proposition \ref{prop:boxx}$ii),$ 
given $t \in \Ia(y^v),$ the $n$-tuple
$(y_1^v t^{h_1}, \ldots, y_n^v t^{h_n})$
corresponds to an element in $\Box(\Sigma^{es}_-).$
Therefore, there exists some $v_k \in I_-$
such that $y^v_k t^{h_k}=1,$ so 
$\gamma^v_k(t) = \gamma^v_k + h_k(\frac{1}{2 \pi i} \log t)
 \in \ZZ.$

We can choose an integer $m << 0$ such that 
$(l_k+mh_k)+ \gamma_k^v(t) \in \ZZ_{\geq 0}.$
As a consequence, the set $\sup (l + mh)$
is a subset of $\Fa \cup (I \setminus v_k).$ 
But the elements of the latter set generate 
a maximal cone in $\Sigma_-^{es}(d),$ so 
$l+mh \in \Sa^{es}_{-} (\gamma^v(t)).$
We have shown that 
$\Sa_{+}^{es} (\gamma^v) \subset
\Sa^{es}_{-} (\gamma^v(t)).$ 

The inverse inclusion is obtained by 
reversing the roles of
$\Sigma_\pm$ and by replacing  $t$ with $t^{-1}.$
The result follows.
\end{proof}

\

For any $l \in \LL,$
we define the integers
$m_{+}(l), m_{-,t}(l),$ $t \in \Ia(y^v),$ 
to be
\begin{equation}\begin{split}\label{eq:mmm}
m_{+}(l)&:= \min \{ m \in \ZZ, mh_j + l_j + \gamma_j^v
\in \ZZ_{\geq 0}, \ 
\text{for some $j,$ $v_j \in I_{+}$} \}, \\
m_{-,t}(l)&:=  \max \{ m \in \ZZ, mh_j + l_j^\prime + \gamma_j^v(t)
\in \ZZ_{\geq 0}, \ 
\text{for some $j,$ $v_j \in I_{-}$} \}.
\end{split}\end{equation}
For any $v \in \Box(\Sigma_+)$ and $t \in \Ia(y^v),$
there exists some $j^\prime, j^{\prime\prime}$ with 
$v_{j^\prime} \in I_+, v_{j^{\prime\prime}} \in I_-$ 
such that $\gamma_{j^\prime}^v, \gamma_{j^{\prime\prime}}^v(t)
\in \ZZ.$ Hence, 
the functions $m_{+}, m_{-,t}: \LL \to \ZZ$
are well defined. Moreover, they are 
are piecewise linear 
with a finite number of linear restrictions. 


\begin{proposition}\label{prop:mm}
Let $l$ be an element of the 
lattice  $\LL$ such that $p(l)$ 
belongs to 
the image of the sets $\Sa^{es}_{+} (\gamma^v),$
$\Sa^{es}_{-} (\gamma^v(t))$
under the projection $p: \LL \to \LL^\prime.$

The element $mh +l \in \LL$ belongs to $\Sa^{es}_{+} (\gamma^v)$
if and only if $m \geq m_+(l),$ and to 
$\Sa^{es}_{-} (\gamma^v(t)),$ $t \in \Ia(y^v),$
if and only if $m \leq m_{-,t}(l).$
\end{proposition}

\begin{proof} If $m < m_+(l),$ then $m h_j + l_j + \gamma^v_j
\notin  \ZZ_{\geq 0},$ for all $v_j \in I_+,$ 
so $I_+$ contains all the elements of $\sup (mh+l).$ But 
the elements of $I_+$ do not generate the rays of a cone
in $\Sigma_+,$ hence $mh + l \not\in \Sa^{es}_{+} (\gamma^v).$
A similar argument shows that, for $t \in \Ia(y^v),$ we have that
$mh +l  \not\in \Sa^{es}_{-} (\gamma^v(t))$
if $m < m_{-,t}(l).$

If $m \geq m_+(l),$ let $i$ be such that $mh_i + l_i
+ \gamma^v_j \in
\ZZ_{\geq 0}.$ Since $p(l)$ is in image of 
$\Sa^{es}_{+} (\gamma^v)$ under the projection 
$p: \LL \to \LL^\prime,$ we conclude that the elements 
$v_j \in \Aa \setminus I$ such that $l_j + 
\gamma^v_j \not\in \ZZ_{\geq 0}$ form a separating 
set $\Fa.$ Hence $\Fa \cup (I \setminus \{ v_i \})$ are the rays 
of a maximal simplex in $\Sigma^{es}_+(d),$
so $mh+ l \in \Sa^{es}_{+} (\gamma^v).$ A 
similar argument shows that, for $t \in \Ia(y^v),$
if $m \leq m_{-,t}(l),$ then 
$mh+ l \in \Sa^{es}_{-} (\gamma^v(t)).$
This ends the proof of the proposition.
\end{proof}

\begin{remark}\label{rem:vuch}
If
$v \in \Box(\Sigma_+)\setminus
\Box(\Sigma_+^{es}),$ any maximal
cone of $\Sigma_+$
containing $\sigma(v)$ is left unchanged by the 
modification. Since for any $l \in \Sa^{es}_{+} (\gamma^v)$
we would have that $\sigma(v) \prec \RR_{\geq 0} \sup (l) \prec
\sigma$ with $\sigma \in \Sigma_+^{es}(d),$ we conclude that
$\Sa^{es}_{+} (\gamma^v)
=\emptyset.$
Moreover, 
by proposition \ref{prop:boxx}$i),$ we also have that 
$v \in \Box(\Sigma_-)\setminus \Box(\Sigma_-^{es})$
and $1 \in \Ia(y^v)$ with $\gamma^v(1)=\gamma^v.$ 
Hence 
$$
\Sa_{+} (\gamma^v)=\Sa_{-} (\gamma^v)
\ \text{and} \ 
\Sa^{es}_{+} (\gamma^v) =\Sa^{es}_{-} (\gamma^v)
=\emptyset.$$
\end{remark} 

\

We are now ready to describe 
the analytic continuation of the analytic 
function $\Xi^v_{+} (z,r) : U_{+} \times B^v \to \CC$
given by (see (\ref{eq:psiv})) 
$$
\Xi_{+}^v (z,r)=
\sum_{l \in \Sa_{+} (\gamma^v)} 
\prod_{j=1}^n \frac{z_j^{l_j+\frac{1}{2\pi i} \log_+ r_j}}
{\Gamma(l_j+\frac{1}{2\pi i} \log_+ r_j+1)}.
$$
As mentioned above, 
in order to simplify our notation, we
choose to denote the $n$ possibly distinct log branches 
with the same symbol $\log_+.$

\begin{remark}\label{rem:chl}
For the given $v \in \Box(\Sigma_+),$ the open domain 
$U_{+} \times B^v$
in $\CC^n \times \CC^n$ has been defined  
in section \ref{sec:gkzk}, and we choose the branches 
$\log_+$ such that 
$$\Sa_{+} (\gamma^v) \subset \Ca_{+}.$$ 
The existence of such a choice follows from corollary 
\ref{cor:choice}. \end{remark}

It will also be useful to consider the analytic function 
$(\Xi^{v}_{+})^{es} (z,r) : U_{+} \times B^v \to \CC$
given by 
\begin{equation}\label{eq:xixi}
(\Xi_{+}^{v})^{es} (z,r)=
\sum_{l \in \Sa^{es}_{+} (\gamma^v)} 
\prod_{j=1}^n \frac{z_j^{l_j+\frac{1}{2\pi i} \log_+ r_j}}
{\Gamma(l_j+\frac{1}{2\pi i} \log_+ r_j+1)}.
\end{equation}

As above, for each $t \in \Ia(y^v),$ we can choose $n$
log branches such that 
$$
\frac{1}{2\pi i}\log_- (y^v_j t^{h_j})= \gamma^v_j(t)=
\gamma^v_j
+ h_j \log t,$$
where $\log t$ is a fixed choice of branch.

The analytic functions
$\Xi^{v(t)}_{-}(z,r), (\Xi^{v(t)}_{-})^{es} (z,r) 
:U_{-} \times B^{v(t)} \to \CC$ are then defined analogously,
for $t \in \Ia(y^v)$ such that there exists a corresponding
twisted sector $v(t) \in \Box(\Sigma_-)$ associated
to $(y_1^v t^{h_1}, \ldots, y_n^v t^{h_n}).$ 
Proposition \ref{prop:boxx} shows that 
such $v(t)$ exists for $t=1 \in \Ia(y^v)$ when
$v \in \Box(\Sigma_{+}) \setminus \Box(\Sigma^{es}_{+}),$
and for all $t \in \Ia(y^v)$ when $v \in \Box(\Sigma^{es}_{+}).$

\

The analytic continuation
will be performed along a path of the form 
$(z(u),y^v)$
starting at a point  
$(z_+, y^v)$ in $U_{+} \times B^v$ 
and ending at a point $(z_-, y^v)$ in
$U_{-} \times B^v.$
The path $z(u)=(z_1(u), \ldots, z_n(u))$
is defined so that, for all $u,$ 
$0 \leq u \leq 1,$ 
$$
\arg z_j(u)= \arg (z_+)_j= \arg (z_-)_j,
$$$$
\log |(z_j(u)|=(1 - u)  \log |(z_+)_j| + u \log |(z_-)_j|
$$
for all $j, 1 \leq j \leq n.$ 

The points $z_\pm \in U_{\pm}$ are chosen 
such that 
the conditions A1)--A3) are satisfied.

A1)
According to corollary \ref{cor:conv1}, 
the domains $U_{\pm} \subset \CC^n$ have the form 
\begin{equation*}
\begin{split}
 \big\{ (z_1, \ldots, z_n) \in \CC^n : &(-\log |z_1|, \ldots, -\log |z_n|) 
\in \Ca_{\pm}^\vee + c_{\pm}, 
\\ &(\arg z_1, \ldots, \arg z_n) \in (-\pi,\pi) \times
\ldots \times (-\pi,\pi) \big\}, 
\end{split}
\end{equation*}
for some appropriate $c_{\pm} \in  \Ca_{\pm}^\vee.$ 
Since the cones $\Ca_{+}^\vee$ and $\Ca_{-}^\vee$ are adjacent 
in the 
secondary fan, and the common facet is included in a 
hyperplane determined by the 
element $h=(h_1, \ldots, h_n) \in \LL,$ 
we can choose the 
points $z_+=((z_+)_1, \ldots, (z_+)_n) \in U_{+}$ and 
$z_-=((z_-)_1, \ldots, (z_-)_n) \in U_{-},$ such that 
$$
\arg (z_+)_j= \arg (z_-)_j, 
-\log |(z_+)_j| + \log |(z_-)_j|= A h_j, \ A > 0,
$$
for all $j, 1 \leq j \leq n.$ 

A2)
The points $z_\pm$ such that 
$(-\log |(z_\pm)_1 |, \ldots, -\log |(z_\pm)_n | \in \Ca_J^\dual + c_J,$
where $J$ is the set of common maximal simplices of $\Sigma_\pm,$
and $\Ca_J=\sum_{\sigma \in J} \Ca_\sigma.$ The element 
$c_J,$ which is located deep inside the cone  
$\Ca^{\dual}_J,$ and the associated open domain 
$U_J \subset \CC^n,$ are provided by the results of 
corollary \ref{cor:conv1}. By convexity, 
for any point $z(u)$ on the analytic continuation path, 
we see that 
$(-\log |z_1(u) |, \ldots, -\log |z_n(u)|) \in \Ca_J^\dual + c_J.$

A3)
The analytic continuation path is chosen so that, 
along the path $z(u),$ 
we have that $- 2 \pi< \arg y(u) < 0,$ where the complex number $y(u)$
is given by
$$
y(u):=e^{i \pi \sum_{j, v_j \in I_-} h_j}\prod_{j=1}^n (z_j(u))^{h_j}.$$
Note that, at least for $ - \pi < \arg z_j(u)=\arg (z_+)_j= \arg (z_-)_j < 0,$
we have that 
$$
\arg y(u) = \sum_{j, v_j \in I_+} h_j \arg z_j(u) + \sum_{j, v_j \in I_-} h_j 
(\pi + \arg z_j(u)) < 0,$$
which shows that it is possible to choose the points $z_\pm$ in 
$U_{\pm}$ such that the argument of $y(u)$ is 
between $-2 \pi$ and $0$ for all $u,$ $0 \leq u \leq 1.$

\begin{theorem}\label{thm:mb} 

\

i) For any $v \in \Box(\Sigma_+),$ 
the function $\Xi^v_{+} (z,r) - (\Xi^{v}_{+})^{es} (z,r)$ 
is analytic in the open domain $U_J \times B^v,$ and 
the open domain $U_J$ contains the sets
$U_{\pm}$ and the path $z(u).$

If $1 \not\in \Ia(y^v),$ 
the function is identically zero in the domain
$U_J \times B^v.$ 

If $1 \in \Ia(y^v),$
then the analytic functions 
$\Xi^v_{+} (z,r) - (\Xi^{v}_{+})^{es} (z,r)$
and  
$\Xi^{v(1)}_{-} (z,r) - (\Xi^{v(1)}_{-})^{es} (z,r)$ 
are equal for all $(z,r) \in U_J \times B^v.$

ii) For any $v \in \Box(\Sigma_+^{es}),$
the analytic continuation
along the path $(z(u),y^v)$
of the germ of the analytic function 
$
(\Xi^{v}_{+})^{es} (z,r) 
$ 
at 
$(z_+, y^v) \in U_{+} \times B^v$
is given by the germ at $(z_-, y^v) \in U_{-} \times B^v $ 
of an analytic function defined as follows. 

If $1 \not\in \Ia(y^v),$ the function can be written as
\begin{equation*}
-\sum_{t \in \Ia(y^v)} \int_{C_t}
T(r,\th) (\Xi^{v(t)}_{-})^{es} (z,r \th^h)\, d\th
+\prod_{j, v_j \in I_+} (1- r^{-1}_j) \,
\phi (z,r).
\end{equation*}

If $1 \in \Ia(y^v),$ the function can be written as
\begin{equation*}
\begin{split}
(\Xi_{-}^{v(1)})^{es} (z,r) -
\sum_{t \in \Ia(y^v)} \int_{C_t} T(r,\th) 
(\Xi^{v(t)}_{-})^{es} &(z,r \th^h)\, d\th  \\
&+ \prod_{j, v_j \in I_+} (1- r^{-1}_j) \,
\varphi (z,r).
\end{split}
\end{equation*}
Here, the integration kernel $T(r,\th)$ is defined by
$$
T(r,\th):=\frac{1}{2 \pi i(\th-1)}
\prod_{j, v_j \in I_-}
\frac{1- r_j^{-1}} 
{1-r_j^{-1}\th^{-h_j}},
$$
$\phi(z,r),$
$\varphi(z,r)$ are analytic functions
on $U_{-}\times B^v,$
and the contours $C_t$ are disjoint
circles in the $\th$-plane
centered at the points 
$t \in \Ia(y^v),$ 
counterclockwise oriented, such that, for $r \in B^v,$
all the poles of $T(r,\th)$ 
are contained inside the discs bounded by the contours 
$C_t.$ 
\end{theorem}

\begin{proof} 
$i)$ The series representation of the function
 $\Xi^v_{+} (z,r) - (\Xi^{v}_{+})^{es} (z,r)$ 
is 
$$
\sum_{l \in \Sa_{+} (\gamma^v) \setminus
\Sa^{es}_{+} (\gamma^v)} \
\prod_{j=1}^n \frac{z_j^{l_j+\frac{1}{2\pi i} \log_+ r_j}}
{\Gamma(l_j+\frac{1}{2\pi i} \log_+ r_j+1)}.
$$
If $1 \not\in \Ia(y^v),$ then by proposition
\ref{prop:sss}$i)$ we have that 
$\Sa_{+} (\gamma^v) \setminus
\Sa^{es}_{+} (\gamma^v)= \emptyset,$ 
so the function is identically zero indeed. 

According to the same result, if $1 \in \Ia(y^v),$ then 
$\Sa_{+} (\gamma^v) \setminus \Sa^{es}_{+} (\gamma^v)=
\Sa_{-} (\gamma^v) \setminus \Sa^{es}_{-} (\gamma^v).$
For any $l$ in these sets, there exists a maximal cone 
$\sigma$ belonging to both $\Sigma_+$ and $\Sigma_-.$ 
We can apply corollary \ref{cor:conv1} for both
fans $\Sigma_\pm$ and $J$ their common subset of maximal cones.
Condition A2) imposed on 
the analytic continuation path $z(u)$ ensures that,
for $0 \leq u \leq 1,$ the point 
$(-\log |z_1(u)|, \ldots, -\log |z_n(u)|) \in \Ca_J^\vee + c_J.$
Note that the branches $\log_{\pm} r_j$ 
are identical for $(r_1, \ldots, r_n) \in B^v.$
It follows that the functions 
$\Xi^v_{\pm} (z,r) - (\Xi^{v}_{\pm})^{es} (z,r)$ 
coincide on $U_J \times B^v,$ with $U_J$ the open domain
in $\CC^n$ provided by corollary \ref{cor:conv1}.

$ii)$ To simplify notation, we set 
$$
\lambda_j:=\frac{1}{2 \pi i} \log_+ r_j.$$

For any $l \in \LL,$ $m_+(l)$ was defined 
as the minimum integer $m$ such that
$mh_j + l_j + \gamma_j^v \in \ZZ_{\geq 0},$ for 
some $j$ with $v_j \in I_+$  
(cf. formula (\ref{eq:mmm})). This shows that  
$m_+(l + mh)= m_+(l) + m$ for any $l \in \LL$ and $m \in \ZZ.$
In particular
$m_+(l - m_+(l) h)= 0,$ and $m_+(l - m h) \not=0$ 
for any other integer $m$ such that $m \not=m_+(l).$

Hence, we can define the piecewise linear 
injection $\iota : \LL^\prime \to \LL$ by 
\begin{equation}\label{eq:iota}
\iota(p(l)):=l - m_+(l) h.
\end{equation}
Note that
$p \circ \iota = Id_{\LL^\prime}.$
We have that $m_+(\iota(p(l)))=0,$ and $\iota(p(l))$ 
is the unique element of the fiber $p^{-1} (p(l))$
with this property. For any $l \in \LL,$ it will be 
convenient to introduce the notations 
$$l^\prime: = \iota(p(l)),$$ 
and
$$
\Sa^\prime:=  \iota(p(\Sa_{+}^{es}(\gamma^v)))= 
\iota(p(\Sa_{-}^{es}(\gamma^v(t)))), t \in \Ia(y^v).
$$
We see that
$$
m_+(l^\prime)=0, \ \text{for any $l^\prime \in \Sa^\prime$}.$$

It is worth noting that, for any $l \in  \Sa_{+}^{es}(\gamma^v),$
we have that
$l^\prime = \iota(p(l))$ belongs to $\Sa_{+}^{es} (\gamma^v).$ 
This is an immediate consequence 
of proposition \ref{prop:mm}, since $m_+(l^\prime)=0.$
Moreover, the choice of the branches $\log_+$ according 
to remark \ref{rem:chl} shows that 
$$
\Sa^\prime \subset  \Sa_{+}^{es} (\gamma^v) \subset \Ca_+.$$

For any $l=mh+ l^\prime \in \LL$ 
and any $\lambda \in \CC^k,$ 
we can write that
\begin{equation}\label{eq:llf}
\begin{split}
&\prod_{j=1}^n \frac{z^{l_j+ \lambda_j}_j} {\Gamma (l_j+\lambda_j+1)}
= \prod_{j=1}^n z_j^{l_j^\prime+ \lambda_j} 
\cdot (-1)^{-\sum_{j, v_j \in I_-} l_j^\prime} \\ 
& \cdot \frac{\prod_{j,v_j \in I_-} (\sin (-\pi \lambda_j)/ \pi) 
\Gamma (-mh_j - l_j^\prime -\lambda_j)}
{\prod_{j,v_j \notin I_-} 
\Gamma (mh_j + l_j^\prime + \lambda_j +1)} 
\big((-1)^{\sum_{j, v_j \in I_-} h_j} \prod_{j=1}^n z_j^{h_j}\big)^m.
\end{split}
\end{equation}

Since $r \in B^v,$ we see that the values 
of the parameters $\lambda_j$ are localized 
around $\gamma^v_j$ 
such that the 
only possible integer value for each
$\lambda_j$ is $\gamma^v_j.$ For $\Vert \lambda - \gamma^v \Vert
< \epsilon,$ 
consider the Mellin-Barnes integral 
(of the type analyzed in lemma \ref{lemma:MB})
$$
I:= \frac{1}{2 \pi i}
\int_{a- i \infty}^{a +i \infty} \ I(s) \ ds,$$
with the integrand $I(s)$ given by  
$$
\frac{\prod_{j,v_j \in I_-} (\sin (-\pi \lambda_j)/ \pi) 
\Gamma (-sh_j - l_j^\prime -\lambda_j) \Gamma(-s) \Gamma(1+s)}
{\prod_{j,v_j \notin I_-} 
\Gamma (sh_j + l_j^\prime + \lambda_j +1)} \, (e^{i\pi}y)^s \, ,
$$
where
$$
y:=e^{i \pi \sum_{j, v_j \in I_-} h_j}\prod_{j=1}^n z_j^{h_j},$$
the path of integration is parallel to the imaginary
axis, and $a$ is a strictly negative real number such that
$ \epsilon < | a | < 1.$ In particular, the contour avoids any poles
of the integrand. The hypotheses of the lemma \ref{lemma:MB}
are satisfied, with $H=2$ and $\beta =0.$
The integral is absolutely convergent and  defines
an analytic function of $y$ for $ -2 \pi < \arg y < 0,$ and 
is equal to the sum 
of the residues at poles on the right of the contour for $|y| < \rho,$
and to the negative of the 
sum of the residues at poles on the left of the contour for
$|y| > \rho,$ where in this case 
$\rho=\prod_{j, v_j \in I} |h_j|^{h_j}.$ 

This special form of the Mellin--Barnes integral has been chosen such that
the residue at any pole $s=m \in \ZZ$ is exactly the last line
of formula (\ref{eq:llf}). The other poles of the integrand $I(s)$ are 
the poles of the product
$\prod_{j,v_j \in I_-} 
\Gamma (-sh_j - l_j^\prime -\lambda_j).$
Set 
$s=m + \theta, m \in \ZZ,$ with 
$$\theta:=\frac{1}{2 \pi i} \log \th,$$
where the branch $\log \th$ is the one chosen above, and we have that
$$
\frac{1}{2\pi i}\log_- (y^v_j t^{h_j})= \gamma^v_j(t)=
\gamma^v_j
+ h_j \log t,$$
for $t \in \Ia(y^v).$
We see that each pole of the 
integrand is the sum of some $m \in \ZZ$ and a complex number
$\theta$ such that
$\lambda_j + h_j \theta \in \ZZ,$ for some $j$ with 
$v_j \in I_-.$ This is equivalent to the condition that 
$$\frac{1}{2 \pi i}(\log_+ r_j  + h_j \log \th)=
\frac{1}{2 \pi i} \log_- (r_j \hat{t}^{h_j}) \in \ZZ,$$
for some $j$ with $v_j \in I_-.$ Let $\Ia(r)$
be the set of such values $\th.$ 

Recall that the set $\Ia(y^v)$ consists of those 
roots of unity $t$ that satisfy 
$\frac{1}{2 \pi i} \log_- (y^v_j t^{h_j}) \in \ZZ,$
for some $j$ with $v_j \in I_-.$ Note that 
for values $r \in (\CC^*)^n$ in an open infinitesimal 
neighborhood of $y^v \in (\CC^*)^n,$ the elements of
the set $\Ia(r)$ are clustered around the elements 
of $\Ia(y^v).$ More precisely, for each $t \in \Ia(y^v),$ 
we can choose mutually disjoint discs 
centered at $t,$ such that any element of the set 
$\Ia(r)$ is contained in one of these discs. For each
$t \in \Ia(y^v),$ the contour $C_t$ is the boundary 
of the corresponding disc. If $1 \in \Ia(y^v),$
we choose the circle $C_1$ centered at $1$ which
contains only the pole $\th=1$ inside. 

The Mellin-Barnes integral $(2 \pi i) I$ introduced above is then equal to
\begin{equation*}\begin{split}
-
\sum_{t \in \{1\}\cup \Ia(y^v)} \sum_{m \geq 0}  
\int_{C_t} &
\frac{\prod_{j,v_j \in I_-} (\sin (-\pi \lambda_j)/ \pi) \,
\Gamma (-h_j(m + \theta) - l_j^\prime -\lambda_j)}
{\prod_{j,v_j \notin I_-} 
\Gamma (h_j(m+\theta) + l_j^\prime + \lambda_j +1)} \\
&\cdot \Gamma(-m-\theta) 
\Gamma(1+m+\theta) (e^{i\pi}y)^{m+ \theta} \, \frac{d\th}{2 \pi i \th},
\end{split}\end{equation*}
when $|y| < \rho, $ and to 
\begin{equation*}\begin{split}
\sum_{t \in \{ 1 \} \cup \Ia(y^v)} \sum_{m < 0}  
\int_{C_t} &
\frac{\prod_{j,v_j \in I_-} (\sin (-\pi \lambda_j)/ \pi) \,
\Gamma (-h_j(m + \theta) - l_j^\prime -\lambda_j)}
{\prod_{j, v_j \notin I^e_-} 
\Gamma (h_j(m+\theta) + l_j^\prime + \lambda_j +1)} \\
&\cdot \Gamma(-m-\theta) 
\Gamma(1+m+\theta) (e^{i\pi}y)^{m+ \theta} \, \frac{d\th}{2 \pi i \th},
\end{split}\end{equation*}
when $|y| > \rho.$ Of course, 
the closed contours $C_t$ avoid all the 
poles of the integrand when $\lambda$ is localized near $\gamma^v.$

In both cases, a direct application of the $\Gamma$--identity, 
shows that, for each $m \in \ZZ,$ the above integrand in $\th$ is equal to
\begin{equation*}\begin{split}
\frac{\pi e^{-i \pi \theta}}{2 \pi i \sin (-\pi \theta)}
\prod_{j, v_j \in I_-} 
&\frac{(-1)^{-h_j m -l_j^\prime}
\sin (-\pi \lambda_j)}{\sin (-\pi(\lambda_j+h_j \theta))}  \\
& \cdot \prod_{j=1}^n \frac{1}
{\Gamma(h_j (m + \theta) + l_j^\prime+ \lambda_j +1)} \, y^{m+\theta}
\end{split}
\end{equation*}
\begin{equation*}\begin{split}
=- 2 \pi i e^{-\pi i \sum_{j, v_j \in I_-} h_j \theta}
&(-1)^{\sum_{j, v_j \in I_-} (-h_j m -l_j^\prime)}  \\
&\cdot T (r, \th)  \prod_{j=1}^n \frac{1}
{\Gamma(h_j (m + \theta) + l_j^\prime+ \lambda_j +1)} \,
y^{m+\theta}
\end{split}\end{equation*}
$$
= - 2 \pi i (-1)^{-\sum_{j, v_j \in I_-} l_j^\prime} 
T (r, \th) \prod_{j=1}^n \frac{z_j^{h_j(m+\theta)}}
{\Gamma(h_j (m + \theta) + l_j^\prime+ \lambda_j +1)},
$$
where, by the statement of the theorem,
the function  
$2 \pi i T(r, \th)$ is equal to
$$ 
\frac{1}{e^{2 \pi i \theta}-1}
\prod_{j, v_j \in I_-}
\frac{1- e^{-2 \pi i \lambda_j}} 
{1-e^{- 2 \pi i (\lambda_j + h_j \theta)}}=
\frac{1}{\th-1}
\prod_{j, v_j \in I_-}
\frac{1- r_j^{-1}} 
{1-r_j^{-1}\th^{-h_j}}.
$$
Therefore, the Mellin-Barnes integral $I$ introduced above is equal to
\begin{equation*}\begin{split}
&(-1)^{-\sum_{j, v_j \in I_-} l_j^\prime}  \\
&\cdot \sum_{t \in \{1\} \cup \Ia(y^v)}
\sum_{m \geq 0} \int_{C_t}  
T (r, \th) \prod_{j=1}^n \frac{z_j^{h_j(m+\theta)}}
{\Gamma(h_j (m + \theta) + l_j^\prime+ \lambda_j +1)} \,d \th,
\end{split}
\end{equation*}
when $|y| < \rho,$ and to 
\begin{equation*}\begin{split}
&
-(-1)^{-\sum_{j, v_j \in I_-} l_j^\prime}  \\
&\cdot \sum_{t \in \{1\} \cup \Ia(y^v)}
\sum_{m < 0} \int_{C_t}  
T (r, \th) \prod_{j=1}^n \frac{z_j^{h_j(m+\theta)}}
{\Gamma(h_j (m + \theta) + l_j^\prime+ \lambda_j +1)} \,d \th,
\end{split}\end{equation*}
when $|y| > \rho.$ It follows that
the analytic continuation of the sum of the former series 
along the path $z(u)$ is the sum of the latter series. 

Recall now that the cones 
$\Ca_{-}^\dual$ and $\Ca_{+}^\dual$ have a common facet 
(maximal dimensional face) which we denote by 
$\tilde{\Ca}.$ It is the unique
facet of the cones $\Ca_{\pm}^\dual$
orthogonal to the element $h \in \LL.$ In order to proceed 
with the analytic continuation, we need two lemmas.

\begin{lemma}\label{lemma:estm}
For any real constants $k,A >0,$ there exists 
an element $\tilde{c}$ deep in the interior of the cone 
$\tilde{\Ca}$ such that, for any $l^{\prime} \in \Sa^{\prime},$
we have that 
$$
\langle u, l^\prime \rangle
\geq k \Vert l^\prime \Vert,
$$
for any 
$u \in \tilde{\Ca} + \tilde{c} + a,$ and any $a \in \LL \otimes \RR$
with $\Vert a \Vert \leq A.$
\end{lemma}

\begin{proof} (of the lemma) 
Note first that 
$$
\langle a, l^\prime \rangle \geq -A \Vert l^\prime \Vert,
$$
since $\Vert a \Vert \leq A.$

The dual of projection  
$p : \LL \to \LL^{\prime}$ is a lattice embedding 
$(\LL^\prime)^\dual \hookrightarrow \LL^\dual$ such that the image 
of $(\LL^\prime)^\dual \otimes \RR$ is the hyperplane in
$\LL^\dual \otimes \RR$ generated by the cone $\tilde{\Ca}.$ 
Hence, with a slight abuse of notation, 
for any $x \in \tilde{\Ca},$ we can write that
$$
\langle x, l^\prime \rangle = \langle
x, p(l^\prime) \rangle,
$$
where we use identical bracket notations for the pairings between
$\LL^\dual$ and 
$\LL,$ and $(\LL^\prime)^\dual$ and $\LL^\prime,$ respectively.

As noted above, the choice of the branches $\log_+$ implies 
that $\Sa^\prime \subset \Ca_+,$ which shows that 
$$p(\Sa^\prime) \subset p(\Ca_+) \subset \tilde{\Ca}^\dual
\subset \LL^\prime \otimes \RR.$$ For any constant
$M >0,$ we can then choose 
an element $\tilde{c}$ deep enough in the interior 
of the cone $\tilde{\Ca}$ (see the end 
of the proof of proposition \ref{prop:cconv})
such that, for any $x$ 
in $\tilde{\Ca}+ \tilde{c}$ we have that
$$
\langle x, p(l^\prime) \rangle \geq M \Vert p(l^\prime) \Vert,$$
for any $l^\prime \in \Sa^\prime.$ But $l^\prime=
\iota(p(l^\prime))$ where the piecewise linear injection 
$\iota : \LL^\prime \to \LL$ has been 
defined by formula (\ref{eq:iota}). Hence, there exists a
constant $K>0$ such that
$$
\Vert p(l^\prime) \Vert \geq K \Vert l^\prime \Vert.$$

We conclude that for an appropriate choice of the element 
$\tilde{c} \in \tilde{\Ca}$ we have that 
$$
\langle x, l^\prime \rangle \geq (k+A) \Vert l^\prime \Vert,$$
for any $x \in \tilde{\Ca}+ \tilde{c}$ and $l^\prime \in \Sa^\prime.$
The lemma follows.
\end{proof}

\begin{lemma}\label{lemma:cconvv}
There exists a value $A >0,$ and an element $\tilde{c} \in \tilde{\Ca},$
such that the set 
$$
V_A:= \{ \tilde{\Ca} + \tilde{c} +a, \ \Vert a \Vert < A \},$$ 
intersects the sets $\Ca_{\pm}^\dual + c_{\pm},$
and, such that the integral 
$$\int_{a+ i \infty}^{a - i \infty}
\sum_{l^\prime \in \Sa^\prime}
\prod_{j=1}^n z_j^{l_j^\prime+ \lambda_j} 
\ (-1)^{-\sum_{j, v_j \in I_-} l_j^{\prime}} \ I(s) \ ds$$
is absolutely convergent, and  defines an analytic function 
of $(z,r)$ in an open domain containing the region $U \times B^v$ 
defined by the restrictions 
\begin{equation*}\begin{split}
U:=&\big\{  z=(z_1, \ldots, z_n) \in \CC^n : 
(-\log |z_1|, \ldots, -\log |z_n|) 
\in V_A, \\& - 2 \pi < \arg y < 0, 
(\arg z_1, \ldots, \arg z_n) \in (-\pi ,\pi) \times
\ldots \times (- \pi ,\pi) \big\}, 
\end{split}
\end{equation*}
where, as before, 
$$
y=e^{i \pi \sum_{j, v_j \in I_-} h_j}\prod_{j=1}^n z_j^{h_j}.$$
\end{lemma}

\begin{proof} (of the lemma) The proof of the lemma is very similar
to that of proposition \ref{prop:cconv} and its corollaries. 
The definition of $V_A,$ as well as the  
the restriction A3) 
imposed on the path $z(u)$ show that $U$ is an  
open set in $\CC^n.$ Moreover, it is clear that by choosing 
a large enough $A >0$ we can ensure that the intersection
of $V_A$ with the sets $\Ca_\pm^\dual + c_\pm$ is non-empty. 

The hypotheses of lemma \ref{lemma:est03} are satisfied. Hence, 
on the line $s=a+it, t \in \RR,$
the absolute value of the integrand is bounded above by a positive
constant multiple of
$$
|y|^{a} e^{-(\pi + \arg y)t} 
(|t|+1)^{R+n/2} e^{-\pi |t|}
\sum_{l^\prime \in \Sa^\prime} (4ek)^{\Vert l^\prime \Vert}
e^{\sum l_j^{\prime} \log |z_j|} \ ,
$$
where $R > 0$ is a positive constant determined by $(r_1, \ldots, r_n).$

As in the proof of proposition \ref{prop:cconv}, the above estimate ignores 
the factors $z_j^{\lambda_j}.$ Lemma \ref{lemma:estm} allows
us to choose $\tilde{c}$ deep inside the cone 
$\tilde{\Ca}$ such that there exists some $\epsilon >0,$ 
for which 
$$
(4ek)^{\Vert l^\prime \Vert}
e^{\sum l_j^{\prime} \log |z_j|} \leq e^{-\epsilon \Vert l^\prime 
\vert},$$
for any $l^\prime \in \Sa^\prime$ and $z \in U.$
The argument restriction
on $y$ insures the absolute convergence of the integral. 
\end{proof} 

Note that the region $U \subset \CC^n$ defined in the previous
lemma imposes no restrictions on $|y|,$ i.e. $U$ contains points
$z=(z_1, \ldots, z_n)$ whose 
associated coordinate $y$ has an absolute value that 
is arbitrarily small or large, as needed. This 
observation allows the analytic continuation procedure between
the regions $|y| < \rho$ and $|y| > \rho$ 
to be performed along the path $z(u).$

Hence, lemma \ref{lemma:cconvv} implies 
that the analytic continuation along the path $z(u)$
from the subdomain $|y| < \rho$ of $U \times B^v$ 
of the series
\begin{equation}\label{eq:ser11}
\frac{1}{2 \pi i}
\sum_{t \in \{1\} \cup \Ia(y^v)}
\sum_{l^\prime \in \Sa^\prime} 
\sum_{m \geq 0} \int_{C_t}  
T (r, \th) \prod_{j=1}^n \frac{z_j^{h_j(m+\theta)+ l_j^\prime+ \lambda_j}}
{\Gamma(h_j (m + \theta) + l_j^\prime+ \lambda_j +1)} \,d \th,
\end{equation}
to the subdomain $|y| > \rho$ of $U \times B^v,$ 
is the series 
\begin{equation}\label{eq:ser12}
-\frac{1}{2 \pi i} 
\sum_{t \in \{1\} \cup \Ia(y^v)}
\sum_{l^\prime \in \Sa^\prime}
\sum_{m < 0} \int_{C_t}  
T (r, \th) \prod_{j=1}^n \frac{z_j^{h_j(m+\theta)
+ l_j^\prime+ \lambda_j}}
{\Gamma(h_j (m + \theta) + l_j^\prime+ \lambda_j +1)} \,d \th.
\end{equation}

For convenience, we introduce the notation
$$P_t(l^\prime,m):=\frac{1}{2 \pi i}
\int_{C_t} 
T (r, \th) \prod_{j=1}^n \frac{z_j^{h_j(m+\theta)
+ l_j^\prime+ \lambda_j}}
{\Gamma(h_j (m + \theta) + l_j^\prime+ \lambda_j +1)} \, d\th.
$$

\begin{lemma}\label{lemma:tran}
For any $t \in \Ia(y^v),$
the series 
\begin{equation*}
S_{-,t}:=\sum_{l^\prime \in \Sa^{\prime}} \
\sum_{0 \leq m \leq \max(0,m_{-,t}(l^\prime))} 
P_t(l^\prime,m)
\end{equation*}
is absolutely convergent for $(z,r) \in U \times B^v,$
and defines an analytic function in the domain $U \times B^v.$
The integers are $m_{-,t}$ have been 
defined by formulae (\ref{eq:mmm}).
\end{lemma}

\begin{proof} (of the lemma) 
For any $t \in \Ia(y^v),$
the series $S_{-,t}$ 
is absolutely
convergent in the subdomain of $U \times B^v$ characterized by 
$|y| < \rho,$ since it is  
a subseries of the 
series (\ref{eq:ser12}) which is absolutely 
convergent in that region. Moreover, 
by proposition \ref{prop:sss}$ii)$ we have that
$p(\Sa^{es}_{+}(\gamma^v))=
p(\Sa^{es}_{-} (\gamma^v(t))),$ 
so, if
$l^\prime \in \Sa^{\prime}$
and $m \leq m_{-,t}(l^\prime),$ 
proposition \ref{prop:mm} shows that 
$l + mh \in \Sa^{es}_{+}(\gamma^v(t)).$
Hence, the series $S_{-,t}$ 
is a subseries of the absolutely
convergent series in the subdomain $|y| < \rho$ 
obtained by integrating over $C_t$ 
the series defining the function 
$\Xi^{v(t)}_{-}(z,r \th^h).$ 
In conclusion, the absolute convergence of
$S_{-,t}$ is independent of the magnitude of
$y,$ therefore it holds everywhere in $U \times B^v.$
\end{proof}

As a direct consequence of lemma \ref{lemma:tran} 
we see that, for any $t \in \Ia(y^v),$ 
the analytic continuation along the path $z(u)$
from the subdomain $|y| < \rho$ of $U \times B^v$ 
of the series
\begin{equation}
\label{eq:ser111}
\sum_{l^\prime \in \Sa^\prime} 
\sum_{m \geq 0} 
P_1(l^\prime,m) + 
\sum_{t \in \Ia(y^v) \setminus \{ 1 \} }
\sum_{l^\prime \in \Sa^\prime}
\sum_{m > \max(0,m_{-,t}(l^\prime))} 
P_t(l^\prime,m),
\end{equation}
to the subdomain $|y| > \rho$ of $U \times B^v,$ 
is the series
\begin{equation}
\label{eq:ser112}
-\sum_{l^\prime \in \Sa^\prime}
\sum_{m < 0} 
P_1(l^\prime,m) -
\sum_{t \in \Ia(y^v) \setminus \{ 1 \}}
\sum_{l^\prime \in \Sa^\prime} 
\sum_{m \leq \max(0,m_{-,t}(l^\prime))} 
P_t(l^\prime,m).
\end{equation}
It is important to remember that, for any 
$l^\prime \in \Sa^\prime,$ we have that
$m_+(l^\prime)=0.$

Note that, 
if $t \in \Ia(y^v) \setminus \{ 1 \}$ 
and $m > m_{-,t}(l^\prime),$ then the function in $\th$ 
\begin{equation}\label{eq:fctt}
T (r, \th) \prod_{j=1}^n \frac{z_j^{h_j(m+\theta)
+ l_j^\prime+ \lambda_j}}
{\Gamma(h_j (m + \theta) + l_j^\prime+ \lambda_j +1)}
\end{equation}
has no poles inside the contour $C_t.$ Hence, 
the second term in the formula (\ref{eq:ser111})
is always zero. Furthermore, 
the second term in formula  
(\ref{eq:ser112}) is equal to 
$$
- \sum_{t \in \Ia(y^v) \setminus \{ 1 \}} \int_{C_t}
T(r,\th) (\Xi^{v(t)}_{-})^{es} (z,r \th^h) \, d\th. 
$$

If $1 \notin \Ia(y^v),$ the function (\ref{eq:fctt})
has only the simple pole $\tilde{t}=1$
inside the contour $C_1,$ for all integers $m.$ 
Hence, in this case, the series (\ref{eq:ser111}) is equal to 
$(\Xi^v_{+})^{es}(z,r).$ 

Moreover, if $1 \notin \Ia(y^v),$
and $m < 0=m_+(l^\prime),$ the residue of 
the function (\ref{eq:fctt}) at $\tilde{t}=1$ 
can be written as a product of 
an analytic function on $\CC^n \times B^v$
and the product
$\prod_{j, v_j \in I_+} (1- r^{-1}_j).$
Hence, in this case, 
the first
term in formula (\ref{eq:ser112}) is of the form
$$
\prod_{j, v_j \in I_+} (1- r^{-1}_j) \,
\phi (z,r),
$$
with $\phi(z,r)$ an analytic function
on $U_{-}\times B^v.$
This proves part $ii)$ of the theorem in the case $1 \not\in
\Ia(y^v).$ 

Let's now assume that $1 \in \Ia(y^v).$
We only have to analyze the first terms in the series
(\ref{eq:ser111}) and (\ref{eq:ser112}). 
Note that, if $m > m_{-,1}(l^\prime),$ then
$$
P_1(l^\prime,m)= 
\prod_{j=1}^n \frac{z_j^{m h_j
+ l_j^\prime+ \lambda_j}}
{\Gamma(m h_j + l_j^\prime+ \lambda_j +1)}. 
$$
It is convenient to introduce the notation
$$
R(l^\prime,m):=
\prod_{j=1}^n \frac{z_j^{m h_j
+ l_j^\prime+ \lambda_j}}
{\Gamma(m h_j + l_j^\prime+ \lambda_j +1)}. 
$$

We proceed with 
a case by case analysis according to signs of the integers
$m_{-,1}(l^\prime).$ 
Namely, we write the set $\Sa^\prime$
as a disjoint union of the subsets 
$\Sa_1$ and $\Sa_2,$ where
\begin{equation*}
\begin{split}
\Sa_1 &:= \{ l^\prime \in \Sa^{\prime} \, : \, 
m_{-,1}(l^\prime) \geq 0 \}, \\
\Sa_2 &:= \{ l^\prime \in \Sa^{\prime} \, : \, 
m_{-,1}(l^\prime) < 0 \}.
\end{split}
\end{equation*}

The terms of the series $\sum_{l^\prime \in \Sa^\prime} 
\sum_{m \geq 0} P_1(l^\prime,m)$ 
(the first half of the series (\ref{eq:ser111}))
coming from $l^\prime \in \Sa_1$ add up to
$$
\sum_{l^\prime \in \Sa_1} \Big(
\sum_{m \geq 0} R(l^\prime,m) -
\sum_{0 \leq m \leq  m_{-,1}(l^\prime)} R(l^\prime,m) + 
\sum_{0 \leq m \leq  m_{-,1}(l^\prime)} P_1(l^\prime,m)
\Big).$$
The terms of the series $\sum_{l^\prime \in \Sa^\prime} 
\sum_{m \geq 0} P_1(l^\prime,m)$ 
coming from $l^\prime \in \Sa_2$ add up to
$$
\sum_{l^\prime \in \Sa_2} 
\sum_{m \geq 0} R(l^\prime,m).$$

The terms of the series $-\sum_{l^\prime \in \Sa^\prime} 
\sum_{m < 0} P_1(l^\prime,m)$ 
(the first half of the series (\ref{eq:ser112}))
coming from $l^\prime \in \Sa_1$ add up to
$$
-\sum_{l^\prime \in \Sa_1} 
\sum_{m < 0} P_1(l^\prime,m).$$
The terms of the series $-\sum_{l^\prime \in \Sa^\prime} 
\sum_{m < 0} P_1(l^\prime,m)$ 
coming from $l^\prime \in \Sa_2$ add up to
$$
-\sum_{l^\prime \in \Sa_2} \Big( 
\sum_{m \leq m_{-,1}(l^\prime)} P_1(l^\prime,m)
+ \sum_{m_{-,1}(l^\prime) < m < 0} R(l^\prime,m)
\Big).$$
Note that, since $m_+(l^\prime)=0$ for all $l^\prime \in
\Sa^\prime,$
the second part of this series is of the form
$$
\prod_{j, v_j \in I_+} (1- r^{-1}_j) \,
\varphi_1(z,r),
$$
where $\varphi_1 (z,r)$ is an analytic function in the 
domain $U_{-} \times B^v.$

When we put everything together, the series 
(\ref{eq:ser111}) is equal to
\begin{equation*}\begin{split}
&\sum_{l^\prime \in \Sa_1 \cup \Sa_2} 
\sum_{m \geq 0} R(l^\prime,m) \\
-&\sum_{l^\prime \in \Sa_1}
\sum_{0 \leq m \leq  m_{-,1}(l^\prime)} R(l^\prime,m) + 
\sum_{l^\prime \in \Sa_1}
\sum_{0 \leq m \leq  m_{-,1}(l^\prime)} P_1(l^\prime,m).
\end{split}
\end{equation*}
It is important to note that lemma \ref{lemma:tran}
implies that the last terms of the previous formula 
define analytic functions  
in $U \times B^v,$ so they do not change under analytic 
continuation.

The series (\ref{eq:ser112}) is equal to
$$
-\sum_{l^\prime \in \Sa_1} 
\sum_{m < 0} P_1(l^\prime,m) 
-\sum_{l^\prime \in \Sa_2} 
\sum_{m \leq m_{-,1}(l^\prime)} P_1(l^\prime,m)+
\prod_{j, v_j \in I_+} (1- r^{-1}_j) \,
\phi_1(z,r).
$$
It follows that
the analytic continuation along the path $z(u)$
from the subdomain $|y| < \rho$ of $U \times B^v$ 
of the series
$$
(\Xi^v_+)^{es}(z,r)=
\sum_{l^\prime \in \Sa_1 \cup \Sa_2} 
\sum_{m \geq 0} R(l^\prime,m)$$ 
to the subdomain $|y| > \rho$ of $U \times B^v,$ 
is the series
$$
\sum_{l^\prime \in \Sa_1}
\sum_{0 \leq m \leq  m_{-,1}(l^\prime)} R(l^\prime,m)
-\sum_{l^\prime \in \Sa_1 \cup \Sa_2} 
\sum_{m \leq m_{-,1}(l^\prime)} P_1(l^\prime,m)$$
$$+
\prod_{j, v_j \in I_+} (1- r^{-1}_j) \,
\varphi_1(z,r)
-\frac{1}{2 \pi i} \sum_{t \in \Ia(y^v) \setminus \{ 1 \}} \int_{C_t}
T(r,\th) (\Xi^{v(t)}_{-})^{es} (z,r \th^h) \, d\th. 
$$
Moreover, for $(z,r) \in U_{-} \times B^v,$
there exists an analytic function $\varphi_2(z,r)$
in $U_{-} \times B^v,$ such that
$$
\sum_{l^\prime \in \Sa_1}
\sum_{0 \leq m \leq  m_{-,1}(l^\prime)} R(l^\prime,m) = 
(\Xi^{v(t)}_-)^{es}(z,r)
+\prod_{j, v_j \in I_+} (1- r^{-1}_j) \,
\varphi_2 (z,r), 
$$
where
$$
\Big(\sum_{l^\prime \in \Sa_1}
\sum_{m < 0} + \sum_{l^\prime \in \Sa_2}
\sum_{m \leq m_{-,1}(l^\prime)} \Big)
R(l^\prime,m)=-
\prod_{j, v_j \in I_+} (1- r^{-1}_j) \,
\varphi_2 (z,r).
$$

The proof of the theorem is then finished if we note that
$$
\sum_{l^\prime \in \Sa_1 \cup \Sa_2} 
\sum_{m \leq m_{-,1}(l^\prime)} P_1(l^\prime,m)=
\int_{C_1} T(r,\th) 
(\Xi^{v(1)}_{-})^{es} (z,r \th^h)\, d\th.
$$
\end{proof}

Let 
$$
B:= \bigcup_{v \in \Box(\Sigma_+)} B^v.$$
Consistent with definition \ref{def:kgamma}, we define 
the functions 
$$
\Xi_+, (\Xi_+)^{es}: U_+ \times B
 \to \CC, 
$$
such that $\Xi_+ (z,r)=\Xi_+^v (z,r),$ $(\Xi_+)^{es} (z,r)=
(\Xi_+^v)^{es} (z,r),$ for all $v \in \Box(\Sigma_+)$ 
and $(z,r) \in U_+ \times B^v.$ We also define
$$
\Xi_-, (\Xi_-)^{es}: U_- \times \bigcup_{v \in \Box(\Sigma_\pm)} B^v
 \to \CC, 
$$
such that $\Xi_- (z,r)=\Xi_-^v (z,r),$ $(\Xi_-)^{es} (z,r)=
(\Xi_-^v)^{es} (z,r),$ for all $v \in \Box(\Sigma_\pm)$ 
and $(z,r) \in U_- \times B^v.$ 
For $v \in \Box(\Sigma_+)
\setminus \Box(\Sigma_-),$ we set 
$\Xi^{v}_{-}=(\Xi^{v}_{-})^{es}=0.$ This choice is 
consistent with the freedom in choosing the 
function $\Xi_-$ discussed in remark \ref{rem:gene}.

Let $\Ia \subset \CC^*$ be 
the set of roots of unity of a large enough order such that 
$\Ia(y^v) \subset \Ia$ for all $v \in \Box(\Sigma_+).$ 
When $t \notin \Ia(y^v),$ we see that  
$
\Xi^{v(t)}_{-}=(\Xi^{v(t)}_{-})^{es}=0.$
With these conventions, the results of the 
theorem can be expressed in a more convenient 
way as follows.

\begin{corollary}\label{cor:mb}

\

i) 
The function $\Xi_{+} (z,r) - (\Xi_{+})^{es} (z,r)$ 
is analytic in the open domain $U_J \times B,$ and 
the open domain $U_J$ contains the sets
$U_{\pm}$ and the path $z(u).$
The analytic functions 
$\Xi_{+} (z,r) - (\Xi_{+})^{es} (z,r)$
and  
$\Xi_{-} (z,r) - (\Xi_{-})^{es} (z,r)$ 
are equal for all $(z,r) \in U_J \times B.$

ii) 
The analytic continuation
along the path $(z(u),y)$
of the germ of the analytic function 
$ 
(\Xi_{+})^{es} (z,r) 
$
at 
$(z_+, y) \in U_{+} \times B$
is given by the germ at $(z_-, y) \in U_{-} \times B$ 
of the analytic function
\begin{equation*}
(\Xi_{-})^{es} (z,r)-
\sum_{t \in \Ia} \int_{C_t} 
T(r,\th) \bigr) 
(\Xi_{-})^{es} (z,r \th^h)\, d\th 
+ \prod_{j, v_j \in I_+} (1- r^{-1}_j) \,
\varphi (z,r).
\end{equation*}
Here,
$\varphi(z,r)$ is an analytic function
on $U_{-}\times B,$ while 
the integration kernel $T(r,\th)$ and the contours 
$C_t$ are defined as in the statement of theorem \ref{thm:mb}. 
\end{corollary} 

We finish this section with a definition and notation which 
will be useful in the next section. 

\begin{definition}\label{def:mb}
The analytic continuation operator along the path $z(u)$
from the domain $U_+ \times B$ to the domain $U_- \times B$ 
is called the {\it Mellin--Barnes operator}. For any analytic 
function $\phi(z,r)$ in $U_+ \times B,$ $MB(\phi) (z,r)$ denotes
its analytic
continuation along the path $z(u)$ to the domain $U_- \times B.$ 
\end{definition}

\section{Toric Birational Maps vs. Analytic Continuation}
\label{sec:vs}

As in the previous section, we consider a modification associated with the 
integral relation $h_1 v_1 + \ldots + h_n v_n= 0$ and the corresponding 
circuit $I= \{ v_j, h_j \not= 0 \} \subset \Aa$ that determines a change of the 
fan $\Sigma_+$ into the fan $\Sigma_-.$ Let $\hat v \in N$ be the vector
$$
\hat v:= \sum_{j, v_j \in I_+} h_j v_j = \sum_{j, v_j \in I_-} (-h_j) v_j.
$$
Let $\hat \Sigma$ be a stacky fan refining the fans $\Sigma_\pm$ 
obtained by 
replacing the cones generated in the fans $\Sigma_\pm$ 
by sets of type $I \setminus v_{\pm},$ with $v_\pm \in I_\pm$ 
respectively, 
with cones generated by sets of the type $\hat{v} \cup I \setminus \{v_-, v_+ \}$
for $v_{\pm} \in I_{\pm}.$ It is important 
to note that the possibly non-primitive 
vector
$\hat v$ is part of the information defining the 
stacky fan $\hat\Sigma.$ 
This definition makes sense even if one of the sets 
$I_-$ or $I_+$ has only one element. In that case, the new vector $\hat{v}$ 
replaces the corresponding generator of a one dimensional cone in $\Sigma_-$
or $\Sigma_+.$ 

We have the following diagram of weighted blowdowns (possibly in codimension
one, if either $|I_-|,$ or $|I_+|$ is equal to one):
\begin{equation}
\begin{split}
\xymatrix{
& \PP_{\hat{\Sigma}} \ar[ld]_{f_+} 
 \ar[rd]^{f_-} \\
\PP_{\Sigma_+} & & \PP_{\Sigma_-}
}
\end{split}
\end{equation}
We will study the properties of the ``Fourier-Mukai'' map
$$FM : K_0 (\PP_{\Sigma_-}, \CC) \to K_0 (\PP_{\Sigma_+}, \CC)$$
defined by
$$FM:=  (f_+)_* (f_-)^*.$$ 
We will use the same notation $R_j, 1 \leq j \leq n,$ for the 
$K$--theory classes induced by the vectors $v_j$ in any
of the toric DM stacks $\PP_{\Sigma_\pm}, \PP_{\hat\Sigma}.$
For details about how the correspondence between 
vectors and $K$--classes works see section 4 in \cite{BH}.
In particular, we may happen that $R_j=1$ when the 
vector $v_j$ does not generate a cone in the corresponding fan. 
We denote by $\Rh$ the $K$--theory class in $\PP_{\hat\Sigma}$
induced by the vector $\hat v.$

Let $\sigma$ be a cone of the fan $\Sigma_-$
generated 
by the vectors $v_j$ with $v_j \in J \subset 
\Aa.$ Assume that 
$\sigma$
is not a subcone of any essential maximal 
cone (see definition \ref{def:ess}). 

In other words, 
any maximal cone of $\Sigma_-$ containing $\sigma$
as a subcone is also a cone of $\Sigma_+$ and 
$\hat\Sigma.$ 
It follows that 
the quotient fans $\Sigma_{\pm}/\sigma,$ 
$\hat{\Sigma}/\sigma$ are all unchanged, and 
by proposition 4.2 of \cite{BCS}, they 
define closed toric substacks in the ambient toric stacks. 
Note that the toric substack induced by the 
cone $\sigma$ and the exceptional toric substack
induced by the vector $\hat v$ in $\PP_{\hat \Sigma}$
have empty intersection. 
Hence, 
the restrictions of the maps 
$f_\pm$
to the closed substack $\PP_{\hat{\Sigma}/\sigma}$
of $\PP_{\hat \Sigma}$
are isomorphisms onto their images in
$\PP_{\Sigma_{\pm}},$ the 
closed substacks $\PP_{{\Sigma_\pm}/\sigma}.$
As a direct consequence, the following 
proposition holds:

\begin{proposition}\label{prop:unessfm}
For any polynomial $\phi (r_1, \ldots, r_n) 
\in \CC[r_1, \ldots, r_n],$ we have that 
$$
FM\big(\prod_{j, v_j \in J}
(1-R_j) \phi(R)\big) =  
\prod_{j, v_j \in J}
(1-R_j) \phi(R).
$$
\end{proposition}

\begin{proof} It is enough to check the statement for $\phi$ an
arbitrary monomial
$\prod_{i=1}^n r_i^{m_i}.$ We can also assume that $J$ contains 
no elements of $I.$ If it did, the elements of  
$J \setminus (J\cap I)$ would generate a cone that is not
a subcone of any essential maximal cone of $\Sigma_-,$
and the statement for $J \setminus (J\cap I)$ would imply the one 
for $J.$ 
Note that, since the elements of $J$ generate a cone
in $\Sigma_-,$ it is not possible for $I$ to be a subset of $J.$

According to theorem 9.1 (the general case) and 
corollary 9.4 (the blowdown of codimension one case) in \cite{BH},
we have that $(f_+)_* (1)=1.$ By the projection formula, 
it is then enough to show that the pull-backs of the classes
$
\prod_{j,v_j\in J}(1-R_j) \prod_{i=1}^n R_i^{m_i}$
from $K_0(\PP_{\Sigma_\pm},\CC)$ coincide in
$K_0(\PP_{\hat\Sigma},\CC).$ Proposition 8.1 in \cite{BH} implies that
these pull-backs are written, with a slight abuse of 
notation, as
$$\prod_{j,v_j\in J}(1-R_j)\prod_{i=1}^n R_i^{m_i} \Rh^{m_\pm}$$
for some integers $m_\pm.$

However, our assumption on the set $J$ implies  
that the vector $\hat{v}$ and the elements of $J$
do not generate a cone in $\hat{\Sigma},$ so
theorem 4.10 in \cite{BH} shows that
$$
\prod_{j,v_j\in J}(1-R_j) (1-\Rh) =0$$
in $K_0(\PP_{\hat\Sigma},\CC).$
We conclude that
$$ 
\prod_{j,v_j\in J}(1-R_j) \Rh^{m_\pm} = 
\prod_{j,v_j\in J}(1-R_j)
$$
for any integers $m_\pm.$ This 
ends the proof of the proposition.
\end{proof}

\begin{proposition}\label{prop:essfm}
For any polynomial $\phi (r_1, \ldots, r_n) 
\in \CC[r_1, \ldots, r_n],$ we have that 
$FM \big( \phi (R)  \big)= \big(FM(\phi) (\Ra)\big) (1),$
where the function $FM(\phi)$ is defined by 
$$
FM(\phi)(r):= \phi(r)- \sum_{t \in \Ia}
\int_{C_t} T(r,\th) \phi(r\th^h) \, d\th.$$
\end{proposition}

The integration kernel $T(r,\th),$ the set of roots
of unity $\Ia \subset \CC^*,$
as well as the contours $C_t,$ are those used in 
the statement of corollary \ref{cor:mb}.

\begin{proof} According to theorem 9.1 (the general case) and 
corollary 9.4 (the blowdown of codimension one case) in \cite{BH},
if $\Rh$ denotes the $K$--theory class determined by the vector
$\hat v$ in $K_0(\PP_{\hat \Sigma},\CC),$
we have the following equality of formal power
series in $\th$
$$
(f_+)_* \big(\frac 1 {1-\Rh^{-1}\th}\big) = 
\frac 1{1-\th} - 
\frac \th{1-\th} \cdot
\prod_{j, v_j \in I_-}\frac {1-R_j^{-1}} {1-R_j^{-1}\th^{-h_j}}.
$$

At this point, it is more convenient to work with 
linear operators on $K$--theory, rather than the 
$K$--theory itself. Let $\Rah : K_0(\PP_{\hat \Sigma}, \CC)
\to  K_0(\PP_{\hat \Sigma}, \CC)$ be the linear 
map (and ring endomorphism) given by multiplication 
with the class $\Rh.$ In order to understand the 
spectrum of $\Rah,$ note first that proposition \ref{prop:boxx}
can be applied for the relation $\hat{v} + \sum_{j, v_j \in I_-} h_j
v_j=0$ inducing the modification from the fan $\Sigma_+$
to the fan $\hat \Sigma.$ Proposition \ref{prop:boxx}
implies then that the maximum ideals of $K_0(\PP_{\hat \Sigma},\CC)$
are among the $(n+1)$-tuples of 
roots unity of the form  
$(t, y^v_1 t^{p_1}, \ldots, y^v_n t^{p_n})$ with $v \in \Box(\Sigma_+),$
$t \in \Ia(y^v),$ and $p_j$ is equal to $h_j,$ for $j$ with $v_j \in I_-,$ and 
zero, otherwise. We conclude that the spectrum of $\Rah$ consists 
of roots of unity contained in $\Ia.$

It follows that, for any integer $k,$ the operator $\Rah^k$ admits the Cauchy 
integral representation (see Appendix \ref{app:2}).
$$
\Rah^k= -\frac{1}{2 \pi i} \sum_{t \in \Ia} 
\int_{C_t} \th^{k-1} (I - \th \, \Rah^{-1})^{-1} \, d\th.$$
Thus, for any polynomial in one variable 
$\psi (r) \in \CC[r],$ we have that
$$
\psi(\Rh)= -\frac{1}{2 \pi i} \sum_{t \in \Ia} 
\int_{C_t} \psi(t)\th^{-1} (I - \th \, \Rah^{-1} )^{-1} (1) \, d\th.$$
The push down formula implies then that
\begin{equation*} 
\begin{split}
(f_+)_* &(\psi(\Rh))= 
-\frac{1}{2 \pi i} \sum_{t \in \Ia}
\int_{C_t}
\frac {\psi(\th)\th^{-1}}{1-\th} \, d\th \\
&+ \frac{1}{2 \pi i} \sum_{t \in \Ia}
\int_{C_t} 
\frac{\psi(\th)}{1-\th}
\prod_{j, v_j \in I_-}(I-\Ra_j^{-1}) (I-\th^{-h_j}\Ra_j^{-1})^{-1} 
(1) \, d\th.
\end{split}
\end{equation*}
The important fact to note is that, in the second line of the 
formula above, the contours $C_t, t \in \Ia,$ enclose all 
the values $\th$ where the operators $I- \th^{-h_j} \Ra_j^{-1}$ 
with $j$ such that $v_j \in I_-,$ are not invertible
on $K_0(\PP_{\Sigma_+},\CC).$

We now analyze the behavior under the pull-back $(f_-)^*$
of a monomial class $ \prod R_j^{m_j}$ in $K_0(\PP_{\Sigma_-},\CC),$
with $m_j$ positive integers. 
Proposition 8.1 in \cite{BH} implies that
\begin{equation*}
\begin{split}
(f_-)^*(\prod_{j=1}^n R_j^{m_j}) 
&= \prod_{j, v_j \notin I_+} R_j^{m_j} \prod_{j, v_j \in I_+} 
(R_j \Rh^{h_j})^{m_j} \\
&=  \Rh^{\sum_{j, v_j \in I} h_j m_j}
(f_+)^*(\prod_{j=1}^n R_j^{m_j}),
\end{split}
\end{equation*}
where we have used the observation that the modifications 
from $\Sigma_-$ and $\Sigma_+$ to $\hat\Sigma$ are induced
by the relations $\hat{v}= \sum_{j, v_j \in I_-} (-h_j) v_j$
and $\hat{v}=\sum_{j, v_j \in I_+} h_j v_j,$ respectively. 
Hence, by the projection formula, we obtain that
\begin{equation*}
\begin{split}
&(f_+)_* (f_-)^*(\prod_{j=1}^n R_j^{m_j}) = \frac{1}{2 \pi i}
\sum_{t \in \Ia}
\prod_{j=1}^n R_j^{m_j}  \cdot \\
&\int_{C_t} 
\frac{\th^{\sum_{j, v_j \in I} h_j m_j}}
{1-\th}\big( -\th^{-1} +
\prod_{j, v_j \in I_-}(I-\Ra_j^{-1}) (I-\th^{-h_j}\Ra_j^{-1})^{-1} 
(1) \big)\, d\th.
\end{split}
\end{equation*}
Since $h_j=0$ for those $j$ with $v_j \notin I,$
we see that 
$$
\prod_{j=1}^n R_j^{m_j} \th^{\sum_{j, v_j \in I} h_j m_j}=
\prod_{j=1}^n (R_j \th^{h_j})^{m_j}.$$
We conclude that, 
for any polynomial $\phi(r_1, \ldots, r_n),$ 
we have that
\begin{equation*}
\begin{split}
&FM(\phi(R_1, \ldots, R_n))=
(f_+)_* (f_-)^*(\phi(R_1, \ldots, R_n)) \\
&= \frac{1}{2 \pi i}
\sum_{t \in \Ia}
\int_{C_t} 
\phi(R_1 \th^{h_1}, \ldots, R_n \th^{h_n})
(\frac{\th^{-1}}{\th -1} - 2 \pi i \, T(\Ra,\th)(1) )\, d\th \\
&=\phi(R_1, \ldots, R_n) - 
\sum_{t \in \Ia}
\int_{C_t}
T(\Ra,\th) \phi(\th^{h_1} \Rah_1, \ldots, \th^{h_n} \Rah_n) (1) \, d\th.
\end{split}
\end{equation*}
This ends the proof of the proposition. 
\end{proof}

\begin{remark}\label{rem:polyn} 
The statements of the previous two propositions describe 
the Fourier--Mukai action on polynomial classes
$\phi(R_1, \ldots, R_n)$ in $K_0(\PP_{\Sigma_-},\CC).$
Both results can be easily extended to the case 
when $\phi(r_1, \ldots, r_n)$ is an analytic function 
in the domain $\cup_{v \in \Box(\Sigma_-)} B^v,$
where $B^v \subset \CC^n$ are disjoint open sets around the 
$n$-tuples of roots of unity $y^v$ (see section \ref{sec:gkzk}
for more details on the choice of the sets $B^v$). It is then enough
to choose a polynomial $\psi(r_1, \ldots, r_n)$ (cf.
Appendix \ref{app:2}) such that
the linear operators  
$\phi(\Ra_1, \ldots, \Ra_n)$ and $\psi(\Ra_1, \ldots, \Ra_n)$ 
coincide on $K_0(\PP_{\Sigma_-},\CC),$ and to 
apply the previous two results for the polynomial $\psi.$
\end{remark}

\begin{theorem}\label{thm:comm}
The following diagram is commutative:
\begin{equation*}
\begin{split}
\xymatrix@C=1in{
(K_0(\PP_{\Sigma_+},\CC))^\vee \ar[r]^{MS_{+}}
 \ar[d]_{FM^\vee}
 & {\mathcal Sol}(U_{+})\ar[d]^{MB}\\
(K_0(\PP_{\Sigma_-},\CC))^\vee \ar[r]^{MS_{-}}
 & {\mathcal Sol}(U_{-})
}
\end{split}
\end{equation*}
\end{theorem}

\begin{proof} For an arbitrary linear function
$f : K_0(\PP_{\Sigma_+},\CC) \to \CC,$ we have that 
$$MB(MS_{+} (f)) = f \big(MB (\Xi_{+})(z,\Ra)(1) \big),
$$
and that
$$
MS_{-} (FM^\vee(f)) = f \big(FM (\Xi_{-})(z,\Ra)(1) \big).
$$
We first write
$
\Xi_{\pm}(z,r)= (\Xi_{\pm} - (\Xi_{\pm})^{es})(z,r) + (\Xi_{\pm})^{es}(z,r).$
In the notation of section \ref{sec:achs}, 
for any $v \in \Box(\Sigma_-),$ the analytic function 
$(\Xi^v_{-} - (\Xi^v_{-})^{es})(z,r)$ on $U_- \times B^v$ 
is the sum of a series made out of terms of the type 
$$\varphi(z,r)=
\prod_{j=1}^n \frac{z_j^{l_j+\frac{1}{2\pi i} \log_- r_j}}
{\Gamma(l_j+\frac{1}{2\pi i} \log_- r_j+1)}.
$$
with $l \in \Sa_{-} (\gamma^v) \setminus  \Sa^{es}_{-} (\gamma^v).$
By proposition \ref{thm:kth} $ii),$ 
$K_0(\PP_{\Sigma_-},\CC)$ is a direct sum 
of Artinian local rings obtained by localizing
at the maximal ideals $(R_1-y_1^v, \ldots, R_n - y_n^v)$
corresponding to all elements
$v \in \Box(\Sigma_-),$
$$
K_0(\PP_{\Sigma_-},\CC) = \bigoplus_{v \in \Box(\Sigma_-)}
(K_0(\PP_{\Sigma_-},\CC))_v.$$
Since $\varphi(z,r)$ vanishes for values of $r$ 
outside of $B^v,$ we see that 
$$
\varphi(z,\Ra) = \varphi(z,\Ra) E_v,$$
where $E_v : K_0(\PP_{\Sigma_-},\CC) \to K_0(\PP_{\Sigma_-},\CC)$
is the projection onto the subspace 
$(K_0(\PP_{\Sigma_-},\CC))_v$ corresponding to $v.$ 

Since $l \in \Sa_{-} (\gamma^v) \setminus  \Sa^{es}_{-} (\gamma^v),$
the elements of $\sup(l)$ generate a cone in $\Sigma_-$
that is not a subcone of any essential maximal cone. In particular,
this shows that $v$ is also in $\Box(\Sigma_+)$ with the same
associated $y^v \in (\CC^*)^n.$

There exists an analytic function $\psi(z,r)$ 
in  $U_- \times B^v$ such that 
$$
\varphi(z,r)=\prod_{j, v_j \in J}(1 - r_j) \psi(z,r),$$
where $J$ consists of those $v_j$ in $\sup(l)$ such that 
$l_j+ 1/(2 \pi i) \log_- y^v_j$ is a negative integer. 
In particular, $y^v_j=1$ for such $j.$ We introduce 
the analytic function $\tilde{\varphi}(z,r)$ on $\CC^n \times B^v$
given by 
$$
\tilde{\varphi}(z,r):= \prod_{j, v_j \in \sigma(v) }(1 - r_j)^{-1} 
\psi(z,r).$$

Proposition \ref{prop:unessfm} and remark \ref{rem:polyn} 
imply that 
$$
FM \big( 
\prod_{j, v_j \in \sup(l)} (I -\Ra_j) \tilde{\varphi}(z,\Ra) (1 )\big)=
\prod_{j, v_j \in \sup(l)} (I -\Ra_j) \tilde{\varphi}(z,\Ra) (1),$$
which means that
$$
FM \big( \varphi(z,\Ra)(1) \big)= \varphi(z,\Ra)(1).$$ We 
conclude that 
$$
FM \big( (\Xi^v_{-} - (\Xi^v_{-})^{es})(z,\Ra) (1) \big)=
(\Xi^v_{-} - (\Xi^v_{-})^{es})(z,\Ra) (1),$$
for all $v \in \Box(\Sigma_-).$ Hence
$$
FM \big( (\Xi_{-} - (\Xi_{-})^{es})(z,\Ra) (1) \big)=
(\Xi_{-} - (\Xi_{-})^{es})(z,\Ra) (1).$$

Corollary \ref{cor:mb} $i)$ shows then that 
$$
FM \big( (\Xi_{-} - (\Xi_{-})^{es})(z,\Ra) (1) \big)=
MB \big( (\Xi_{+} - (\Xi_{+})^{es})(z,\Ra) (1) \big).
$$
Moreover, corollary \ref{cor:mb} $ii)$ and proposition
\ref{prop:essfm} (combined with remark \ref{rem:polyn})
show that
$$
FM \big( (\Xi_{-})^{es}(z,\Ra) (1) \big)=
MB \big( (\Xi_{+})^{es}(z,\Ra) (1) \big),
$$
which ends the proof of the theorem. 
\end{proof}

\section{Appendices}

\subsection{Auxiliary Analytic Results and Estimates}
\label{app:1}

In what follows,
we define the norm of $r =(r_1, \ldots, r_n) \in \CC^n$ to be 
the positive real number
$\Vert r \Vert := \vert r_1 \vert + \ldots + \vert r_n \vert.$
As a convention, the arguments of all the complex numbers 
used here are chosen in $[-\pi,\pi].$

\begin{lemma}\label{lemma:jensen}
Let $a_1, \ldots, a_p$ be strictly positive real numbers. 
Then
$$
\prod_{j=1}^p a_j^{-a_j}  \leq 
\big(\frac{\Vert a \Vert}{p}\big)^{-\Vert a \Vert}.
$$
\end{lemma}

\begin{proof} Note that the function $f(x)=-x \log x$
is concave down on $(0, + \infty).$ This means that 
$$
\frac{\sum_{j=1}^n -x_j \log x_j}{n} \leq -\frac{\Vert x \Vert}{n}
\log \frac{\Vert x \Vert}{n}.
$$
After applying the exponential to the two sides of 
the above inequality, the desired result is obtained.
\end{proof}

\begin{lemma}\label{lemma:est01} 

There exists a positive constant $M>0$ such that
$$
| \frac{1}{\Gamma(z)}| \leq M \cdot (|x|+|y|)^{-x+1/2} e^{x+y\theta},
$$
for any complex number $z=x+iy=Re^{i \theta}.$

\end{lemma}

\begin{proof} 
According to Stirling's formula, we have that, for 
a fixed complex number $u$ and for any $\delta > 0,$
\begin{equation}\label{ster2}   
\Gamma(z+u) = (2\pi)^{1/2} \, z^{z+u- 1/2} e^{-z} O(1)
\end{equation}
where, as $|z| \to \infty,$ 
$O(1)$ goes to $1$ {\it uniformly} in
$$
| \arg z | < \pi - \delta.
$$
As a direct consequence, we have that, if we write $z=x+iy=Re^{i\theta},$ 
$|\theta| < \pi - \delta,$ then ($u=0$)
$$
| \frac{1}{\Gamma(z)}|=  (2\pi)^{-1/2} R^{-x+1/2} e^{x+y\theta}
O(1), \ \text{when $R \to \infty$}.
$$
If $- \pi \leq \theta \leq - \pi + \delta,$ or $\pi - \delta \leq \theta \leq
\pi,$ we use the Gamma identity $\Gamma(z)\Gamma(1-z)=\pi/\sin(\pi z),$
$z \notin \ZZ_{\leq 0},$ and again  
Stirling's formula for $\Gamma(1-z)$ with $u=1$ to write that
$$
| \frac{1}{\Gamma(z)}|= e^{\pi |y|}
|e^{-2\pi (|y| \pm ix)} -  1| \, 
 (2\pi)^{-1/2} R^{-x+1/2} e^{x+y \theta^*}  O(1),
$$
where $\theta^*=\arg (-z).$ 
Note, however, that, for any $z=x+iy=Re^{i \theta},$  
$$
\pi |y| + y \theta^*=y \theta.
$$

Clearly, there exist a positive constant 
$a > 0,$ such that
$$
|e^{-2\pi (|y| \pm ix)} -  1| < b,
$$
for any complex number $z=x+iy.$ 

Since 
$$
\frac{1}{\sqrt{2}}(|x|+|y|) \leq R=\sqrt{x^2+y^2} \leq |x|+|y|,
$$
we can replace $R$ by $|x|+|y|$ above, 
and conclude that there exist a positive constant $M >0,$ such that 
$$
|\frac{1}{\Gamma(z)}| \leq 
M \cdot (|x|+|y|)^{-x+1/2} e^{x+y\theta},
$$
for any complex number $z=x+iy=Re^{i\theta}.$ 
This ends the proof of the lemma. 
\end{proof}

\begin{lemma}\label{claim} 
For any $\delta >0,$ 
there exists a positive constant $B> 0,$ such that
$$
\prod_{j=1}^n (|x_j|+|y_j|)^{-x_j+1/2} \leq B \cdot (4n)^{\Vert x \Vert},
$$
for any $x=(x_1, \ldots, x_n), y=(y_1, \ldots, y_n)$ in $\CC^n$  
such that
$$
|x_1 + \ldots + x_n | \leq \delta, \ \Vert y \Vert \leq \delta. 
$$
\end{lemma}

\begin{proof}
Note that
\begin{equation}\label{eq:f1}
\prod_{j=1}^n (|x_j|+|y_j|)^{-x_j+1/2} = 
\prod_{j=1}^n 
\big(\frac{|x_j|+|y_j|}{\Vert x \Vert + \Vert y \Vert}\big)^{-x_j+1/2}
(\Vert x \Vert + \Vert y \Vert)^{\sum (-x_j+1/2)}.
\end{equation}
We see that
$$
\prod_{j=1}^n 
\big(\frac{|x_j|+|y_j|}{\Vert x \Vert + \Vert y \Vert}\big)^{1/2} \leq 1.
$$
If $x_j \leq 0,$ then
$$
\big(
\frac{|x_j|+|y_j|}{\Vert x \Vert + \Vert y \Vert}\big)^{-x_j} 
\leq 1,
$$
while if $x_j > 0,$ then
$$
\big(
\frac{|x_j|+|y_j|}{\Vert x \Vert + \Vert y \Vert}\big)^{-x_j} 
\leq \big(
\frac{x_j}{\Vert x \Vert + \delta}\big)^{-x_j}. 
$$

Let $p$ be the number of strictly positive $x_j,$ and assume 
that $p >0.$
We apply lemma \ref{lemma:jensen} to the positive real numbers
$$
\frac{x_j}{\Vert x \Vert + \delta}, \ x_j >0,
$$
and get that
$$
\prod_{x_j >0} \big(\frac{x_j}{\Vert x \Vert + \delta}\big)^{-x_j} 
\leq \Big(\frac{\sum_{x_j >0} x_j} {p(\Vert x \Vert + \delta)} \Big)
^{-\sum_{x_j >0} x_j}.
$$
Since
$$
\sum_{x_j >0} x_j + \sum_{x_j \leq 0} x_j >- \delta, \ 
\text{and} \ \Vert x \Vert = 
\sum_{x_j >0} x_j - \sum_{x_j \leq 0} x_j,
$$
we see that
$$
\Vert x \Vert < 
2 \sum_{x_j >0} x_j + \delta.
$$
Hence
\begin{equation} \nonumber
\begin{split}
&\Big(\frac{\sum_{x_j >0} x_j} {p(\Vert x \Vert + \delta)} \Big)
^{-\sum_{x_j >0} x_j} \leq
\Big(\frac{\sum_{x_j >0} x_j} {2p( \sum_{x_j >0} x_j + \delta)} \Big)
^{-\sum_{x_j >0} x_j} \leq \\
& \leq (2p)^{\sum_{x_j >0} x_j} (1+ \frac{\delta}{\sum_{x_j>0} x_j})^{\sum_{x_j>0} x_j} \leq 
K \cdot (2n)^{\Vert x \Vert},
\end{split}
\end{equation}
for some positive constant $K>0.$

We conclude that
\begin{equation}\label{eq:f2}
\prod_{j=1}^n
\big(\frac{|x_j|+|y_j|}{\Vert x \Vert + \Vert y \Vert}\big)^{-x_j+1/2}
\leq K \cdot (2n)^{\Vert x \Vert},
\end{equation}
and the inequality obviously holds also in the case when 
there are no positive $x_j$'s. i.e. when $p=0.$

Moreover
$$
(\Vert x \Vert + \Vert y \Vert)^{\sum (-x_j+1/2)} \leq 
(\Vert x \Vert + \Vert y \Vert)^{\pm \delta + k/2}. $$
But $ \Vert y \Vert \leq \delta,$ so there exists a positive
constant $B >0$ (depending on $\delta$) such that
\begin{equation}\label{eq:f3}
(\Vert x \Vert + \Vert y \Vert)^{\sum (-x_j+1/2)} \leq 
B \cdot 2^{\Vert x \Vert}, 
\end{equation}
for any value of $\Vert x \Vert \geq 0.$

By combining the formulae (\ref{eq:f1}), (\ref{eq:f2}) and (\ref{eq:f3}), 
we can write
$$
\prod_{j=1}^n (|x_j|+|y_j|)^{-x_j+1/2} \leq K \cdot B \cdot (4n)^{\Vert x \Vert},
$$
which proves the lemma. 
\end{proof}

\begin{lemma}\label{lemma:est02} 
For any $\delta >0,$ 
there exists a positive constant $A> 0,$ such that
$$\big\vert
\prod_{j=1}^n \frac{1} {\Gamma (x_j+iy_j)} \big\vert
\leq
A \cdot (4n)^{\Vert x \Vert},
$$
for any $x=(x_1, \ldots, x_n), y=(y_1, \ldots, y_n)$ in $\CC^n$  
such that
$$
|x_1 + \ldots + x_n | \leq \delta, \ \Vert y \Vert \leq \delta. 
$$
\end{lemma}

\begin{proof}

Clearly, we can safely assume that $\Vert x \Vert + \Vert y \Vert >0.$
According to lemma \ref{lemma:est01}, we have that
$$\big\vert
\prod_{j=1}^n \frac{1} {\Gamma (x_j+iy_j)} \big\vert
\leq M^n \ 
e^{\sum x_j + \sum y_j \theta_j} 
\prod_{j=1}^k (|x_j|+|y_j|)^{-x_j+1/2}.
$$
Note first that 
$$
\sum x_j + \sum y_j \theta_j \leq \delta + \pi \delta.
$$

According to the previous lemma, 
there exists a constant $B > 0$ such that
$$
\prod_{j=1}^n (|x_j|+|y_j|)^{-x_j+1/2} \leq B \cdot (4n)^{\Vert x \Vert}.
$$

Hence 
$$
\big\vert
\prod_{j=1}^n \frac{1} {\Gamma (x_j+iy_j)} \big\vert
\leq M^n \, e^{\delta + \pi \delta} \, B \cdot (4n)^{\Vert x \Vert},
$$
which ends the proof of the lemma.
\end{proof}

\begin{lemma}\label{lemma:est03} 
Let $h=(h_1, \ldots, h_n)$ 
be a fixed element in $\RR^n,$ and $I_-, I_+$ two subsets 
that determine a partition of  
$\{ 1, \ldots, n \}.$ Define
$$
H:=\sum_{j \in I_-} |h_j| - \sum_{j \in I_+} |h_j|.
$$ 

For any $\epsilon >0$ and $\delta >0,$ 
there exists a positive constant $A> 0,$ such that
$$\big\vert 
\frac{\prod_{j \in I_-} \Gamma (1-x_j-iy_j - i(h_j t))}
{\prod_{j \in I_+} \Gamma (x_j+iy_j + i(h_j t))}
 \big\vert
\leq
A \cdot (|t|+1)^{\delta+n/2} (4en)^{\Vert x \Vert}  e^{-\pi H |t|/2},
$$
for any $t \in \RR,$ and 
any $x=(x_1, \ldots, x_n), y=(y_1, \ldots, y_n)$ in $\CC^n$  
such that
$$
|x_1 + \ldots + x_n | \leq \delta, \ \Vert y \Vert \leq \delta, 
$$
and such that, for all $j \in I_-,$ 
the complex numbers $1-x_j-iy_j - i(h_j t)$
are located at a distance 
greater than $\epsilon$ from any integer.
\end{lemma}

\begin{proof}
For $j \in I_-,$
we use the Gamma identity 
to write 
$$
|\Gamma (1-x_j-iy_j - i(h_j t))|= 
\frac{2 \pi e^{-\pi |y_j+h_j t|}} {|e^{-2\pi (|y_j +h_j t| \pm ix_j)} -  1|} \ 
\big| \frac{1}{\Gamma (x_j+iy_j + i(h_j t))} \big|.
$$

The hypothesis that $1-x_j-iy_j - i(h_j t)$ are located at a distance 
greater than $\epsilon$ from any integer, for
all $j \in I_-,$ guarantees that 
$$1/|e^{-2\pi (|y_j +h_j t| \pm ix_j)}-1|$$
is bounded from above for all $j \in I_-.$

If we restrict the real number $t$ to 
a bounded range $|t| \leq \Lambda,$ then lemma \ref{lemma:est02} 
provides the required result with a constant $A$ that {\it depends} on
$\Lambda.$ It remains to understand what happens for $|t| > \Lambda,$
for some fixed $\Lambda >1.$

Note that
$$
-|y_j+h_j t| \leq |y_j| - |h_j t|.
$$
As a result of the above considerations,
the quotient of products of Gamma functions 
invoked in the statement of the lemma is bounded above by
a constant multiplied by 
$$
e^{-\pi \sum_{j \in I_-} |h_j t|} \prod_{j=1}^n
\big| \frac{1}{\Gamma (x_j+iy_j + i(h_j t))} \big|.
$$
According to lemma \ref{lemma:est01}, this expression is 
bounded by a constant multiplied by 
$$
e^{-\pi \sum_{j \in I_-} |h_j t| +\sum_{j=1}^n (x_j+(y_j+h_j t)\theta_j)}
\prod_{j=1}^n (|x_j|+|y_j|+|h_j t|)^{-x_j+1/2},
$$
with $x_j+iy_j+ ih_j t=R_j e^{i \theta_j}.$ Since it is the 
case that $(y_j+h_j t)\theta_j=|y_j + h_j t| |\theta_j| 
\leq (|y_j| + |h_j t|) |\theta_j| \leq \delta \pi + 
|h_j t| |\theta_j|,$ and $\sum_{j=1}^n x_j < \delta,$  
we infer that the above expression 
is also bounded above by a constant multiplied by 
\begin{equation}\label{eq:form1}
e^{(-\pi \sum_{j \in I_-} |h_j| + \sum_{j=1}^n |h_j \theta_j|)|t|}  
\prod_{j=1}^n (|x_j|+|y_j|+|h_j t|)^{-x_j+1/2}.
\end{equation}

We first study the sum  $\sum_{j=1}^n |h_j \theta_j| |t|.$
Note that 
for any $x,s$ with $s \not= 0$ one has
$$
|s| \arctan(|x|/|s|) \leq |x|,
$$
which shows that 
$$
(|arg( x + is)| - \pi/2) |s| \leq |x|
$$
for all $x$ and $s$.
When we apply this to $x=x_j$ and $s=y_j+h_jt$,
we get
$$
|\theta_j||y_j+h_jt| \leq \pi/2 |y_j+h_j t| + |x_j|
$$
which implies

$$
|\theta_j||h_j||t| \leq \pi/2 |h_j||t| + |x_j| + C,
$$
where the constant $C$ depends on $\delta$ only.
We sum this over all $j$ and exponentiate to get that
the factor
$$e^{(-\pi \sum_{j \in I_-} |h_j| + \sum_{j=1}^n |h_j \theta_j|)|t|}$$ 
in formula (\ref{eq:form1})
is bounded
above by a constant (independent of the $x_j$'s) multiplied by 
$$
e^{-\pi H |t|/2} e^{\Vert x \Vert}.
$$
where $H=\sum_{j \in I_-} |h_j| - \sum_{j \in I_+} |h_j|.$

Let's now analyze the other factor of the formula (\ref{eq:form1}) for $|t| > \Lambda.$
We see that
\begin{equation*}
\begin{split}
&\prod_{j=1}^n (|x_j|+|y_j|+|h_j t|)^{-x_j+1/2}= \\
=&|t|^{\sum (-x_j+1/2)} 
\prod_{j=1}^n (|x_j/t|+|y_j/t|+|h_j|)^{-x_j+1/2}.
\end{split}
\end{equation*}
Since $|t| > \Lambda > 1,$ we see immediately that
$$
|t|^{\sum (-x_j+1/2)} \leq (|t|+1)^{\delta+n/2}.$$ 

According to lemma \ref{claim}, we have that 
$$
\prod_{j=1}^n (|x_j/t|+|y_j/t|+|h_j|)^{-x_j/t+1/2} \leq
B \cdot (4n)^{\Vert x/t \Vert} \leq B \cdot (4n)^{\Vert x\Vert /\Lambda}. 
$$

Hence the factor $\prod_{j=1}^n (|x_j|+|y_j|+|h_j t|)^{-x_j+1/2}$ is bounded 
above by a constant multiplied by 
$$
(|t|+1)^{\delta+n/2} (4n)^{\Vert x\Vert},
$$
which ends the proof of the lemma.

\end{proof}

The following property is essentially stated in \cite{bateman}, page 49,
and \cite{WW}, \S 14.5. We include a proof for completeness. 

\begin{lemma}\label{lemma:MB}
Consider the integral
\begin{equation}\nonumber
\int_{\gamma+ i \infty}^{\gamma -i \infty}
 \frac{\prod \Gamma(A_j s + a_j) \prod \Gamma(-C_j s + c_j)}
{\prod \Gamma(B_j s + b_j) \prod \Gamma(-D_j s + d_j)} \; y^s \; ds,
\end{equation}
with $\gamma$ real, and $A_j, B_j, C_, D_j$ all real and strictly positive.
The path of integration
is parallel to the imaginary axis for large $|s|,$ but it
can be curved elsewhere so that it avoids the poles of the
integrand. Introduce the following notations
\begin{equation}
\begin{split}\nonumber
H:&=\sum A_j  + \sum C_j - \sum B_j- \sum  D_j,\\
\beta :&= \sum A_j  - \sum C_j - \sum B_j + \sum  D_j,
\\
\eta  :&= \Re \big(\sum (a_j  - \frac{1}{2}) + \sum (c_j  - \frac{1}{2})
- \sum (b_j  - \frac{1}{2}) - \sum (d_j  - \frac{1}{2}) \big),\\
\rho :&= (\prod A_j ^{-A_j}) (\prod C_j ^{C_j})
(\prod B_j ^{B_j}) (\prod D_j ^{-D_j}).
\end{split}
\end{equation}

i) For 
\begin{equation}\nonumber
s=\gamma + it, \ y=R e^{i\theta},
\end{equation}
the absolute value of the integrand has the asymptotic form
\begin{equation*}
e^{-\frac{1}{2} H \pi \vert t \vert } \; \vert t \vert\
^{\beta \gamma +\eta}\; R^\gamma\; e^{-\theta t}\; \rho^{-\gamma}
\end{equation*}
when $|t|$ is large.
Therefore, if
\begin{equation}\nonumber
H > 0,
\end{equation}
the integral is absolutely convergent (and defines
an analytic function of $y$) in any domain contained in
\begin{equation}\nonumber
\vert \arg y \vert < \min(\pi, \frac{H\pi}{2}).
\end{equation}

ii) 
Moreover, if $\beta=0,$ the integral is equal to the sum
of the residues on the right of the contour for $|y| < \rho,$
and to the negative of the 
sum of the residues  on the left of the contour for
$|y| > \rho$ (these facts are obtained by closing the contour to
the right, respectively to the left, with a semicircle of
radius $r \to \infty$).
\end{lemma}

\begin{proof} Part $i)$ is a direct consequence of Stirling's formula (\ref{ster2}).

For part $ii),$ 
assume that $|y| < \rho.$ It is enough to show that the integral
over a semicircle $C$ of radius $M > 0$ on the right of the 
imaginary axis and centered at the origin goes to zero when 
$M$ goes to $+ \infty.$ It is possible to choose a sequence 
$M_n \to + \infty$
while making sure that the expressions $-C_j s + c_j$ 
and  $-D_j s + d_j$ for $s=M_n e^{i\theta}, -\pi/2 \leq \theta \leq \pi/2,$ 
are at a distance greater than some $\epsilon>0$ from all negative
integers. 

The integrand can be written as 
$$ 
I(s):= \frac{\prod \Gamma(A_j s + a_j) \prod \Gamma(D_j s +1- d_j)}
{\prod \Gamma(B_j s + b_j) \prod \Gamma(C_j s + 1-c_j)} \, 
\frac{\prod \sin(\pi (-D_j s +d_j))}{\prod \sin(\pi (-C_j s +c_j))}
\, y^s,
$$
Note that we're working under the assumptions that $H>0$
and $\beta=0.$ It follows that
$$
\sum C_j - \sum D_j = \sum A_j -  \sum B_j =H/2 >0.
$$
      
Stirling's formula (\ref{ster2}) implies that, for $H=0$ 
and $s=M_n e^{i\theta}$
$$
O (M_n^{\eta}) \, \frac{\prod \sin(\pi (-D_j s +d_j))}
{\prod \sin(\pi (-C_j s +c_j))} \; ({y/\rho})^s,
$$
where the symbol $O$ is independent of $\theta=\arg s$ when 
$s$ is on the semicircle.

Note that we write $a_n=O(b_n)$ for two sequences $(a_n), (b_n),$ if 
$$
| a_n /b_n | < K, n >>0,
$$ 
with $K$ independent of $n.$
     
We have that
$$
4| \sin (x +iy) |^2 = e^{2y} + e^{-2y} -2 \cos (2x), \ x,y \in {\mathbb R}.
$$
The choice explained above of the semicircles $s=M_n e^{i\theta}, 
-\pi/2 \leq \theta \leq \pi/2,$ guarantees that 
$$
\vert \frac{1}{\prod \sin(\pi (-C_j s +c_j))} \vert < K, 
$$
with $K$ independent of $n.$

It follows that
$$
\frac{\prod \sin(\pi (-D_j s +d_j))}
{\prod \sin(\pi (-C_j s +c_j))} =
O \big( \exp \big( (\sum D_j-\sum C_j) M_n \pi |\sin \theta| \big) \big),
$$
and, for $| y/\rho | < 1,$
$$
( y/ \rho)^s =O \big( \exp \big( 
{ M_n \cos \theta \log |y/ \rho| - M_n \sin \theta
\arg y} \big) \big).
$$
But $\sum D_j-\sum C_j= -H/2,$ so
$$
\frac{\prod \Gamma(A_j s + a_j) \prod \Gamma(D_j s +1- d_j)}
{\prod \Gamma(B_j s + b_j) \prod \Gamma(C_j s + 1-c_j)} \, 
\frac{\prod \sin(\pi (-D_j s +d_j))}{\prod \sin(\pi (-C_j s +c_j))}
\, y^s = 
$$
$$
= O \big( 
M_n^{\eta} \, \exp \big( M_n(-\frac{1}{2}H \pi  |\sin \theta|
- \arg y \sin \theta) 
+M_n \cos \theta \log |y/ \rho|  \big) \big). 
$$
Choose $\delta >0$ such that $\delta < \frac{1}{2}H \pi \pm 
\arg y,$
and we see that the above expression is in fact 
$$
O \big( 
M_n^{\eta} \, \exp \big( - \delta M_n  |\sin \theta|
+M_n \cos \theta \log |y/ \rho|  \big) \big). 
$$
Hence, for $| y/ \rho| < 1,$ and $-\pi/2 \leq \theta \leq 
\pi/2,$ the integrand tends to zero
sufficiently rapidly (when $n \to \infty$) to ensure that
the integral along the semicircle tends to zero.
\end{proof}

\subsection{Functions of linear operators}
\label{app:2}

We briefly recall some facts from the spectral theory of linear operators
on finite dimensional vector spaces. The details can be found in section VII.1
of \cite{DS}, where the case of one linear operator is treated. In our discussion, 
$\Ra_1,\ldots, \Ra_n$ are mutually commuting linear operators on a finite dimensional 
complex vector space $V,$ i.e
$\Ra_i \Ra_j=\Ra_j \Ra_i$ for all $i,j,$ $1 \leq i,j \leq n.$
If 
$P(x_1, \ldots, x_n)$ is a polynomial with complex coefficients,
then $P(\Ra_1, \ldots, \Ra_n)$ is a well defined linear operator on $V.$ For a
linear operator $\Ra$ on $V,$ 
its 
{\it spectrum} $s(\Ra)$ is the set of complex numbers $\lambda$ such that 
$\Ra- \lambda I $ is not one-to-one. The {\it index} $\nu(\lambda)$ of a complex
number $\lambda$ is the smallest non-negative integer $\nu$ such that 
$$
\{ x \, \vert \,(\Ra - \lambda I)^\nu x=0 \} = 
\{ x \, \vert \, (\Ra - \lambda I)^{\nu+1} x=0 \}.$$ 
Given two complex polynomials $P,Q,$ in $n$ variables
then 
$$P(\Ra_1,\ldots,\Ra_n)=Q(\Ra_1,\ldots,\Ra_n)$$ 
if and only if 
$P-Q$ is divisible by $(x_1-\lambda_1)^{\nu(\lambda_1)} \ldots (x_n-\lambda_n)^{\nu(\lambda_n)}$
for any $(\lambda_1, \ldots, \lambda_n) \in 
s(\Ra_1) \times \ldots \times s(\Ra_n) \subset \CC^n.$
This property provides 
the definition of the linear operator $f(\Ra_1, \ldots,\Ra_n)$ for every function 
$f : \CC^n \to \CC$ which
is analytic in an open domain (not necessarily connected!) 
that contains $s(\Ra_1) \times \ldots \times s(\Ra_n)
\subset \CC^n.$ Indeed, it is enough to consider
a polynomial $P(x_1, \ldots, x_n)$ such that 
$$
\partial_1^{m_1} \ldots \partial_n^{m_n}
P (\lambda_1, \ldots, \lambda_n)= 
\partial_1^{m_1} \ldots \partial_n^{m_n}
f (\lambda_1, \ldots, \lambda_n),$$
for all $\lambda_j \in s(\Ra_j)$
and $0 \leq m_j \leq \nu(\lambda_j)-1,$ and set 
$$f(\Ra_1, \ldots, \Ra_n):= P(\Ra_1, \ldots, \Ra_n).$$
It follows that, for such a function $f,$ the linear operator 
$f(\Ra_1, \ldots, \Ra_n)$ can be expressed as 
$$\sum_{\lambda_1\in s(\Ra_1),\ldots,
\lambda_n \in s(\Ra_n)} \ \
\sum_{0\leq j_1 < \nu(\lambda_1), \ldots ,0 \leq j_n < \nu(\lambda_n)} $$
$$\frac{(\Ra_1-\lambda_1I)^{j_1}E_1(\lambda_1)}{j_1 !} \ldots 
 \frac{(\Ra_n-\lambda_n I)^{j_n}E_n (\lambda_n)}{j_n !}\,
\partial_1^{j_1} \ldots \partial_n^{j_n}
f (\lambda_1, \ldots, \lambda_n),
$$ 
where the operators $E_j (\lambda)$ are the usual projections onto 
the kernels of the operators $(\Ra_j-\lambda I)^{\nu(\lambda)}$
for some eigenvalue $\lambda \in s(\Ra_j)$ (compare to 
theorem VII.1.8 in \cite{DS}). As a direct consequence, there
exists a Cauchy type integral representation of the operator
$f(\Ra_1, \ldots, \Ra_n).$ Namely, assume that $f$ is 
analytic in an open domain in $\CC^n$ containing 
$U_1 \times \ldots \times U_n,$ where, for all $j, 1 \leq j \leq n,$
$U_j$ contains $s(\Ra_j),$ and 
the (positive oriented) boundary $B_j$ of $U_j$ consists of a 
finite union of closed rectifiable Jordan curves. Then
$
f(\Ra_1, \ldots, \Ra_n)$ is given by 
$$\frac{1}{(2 \pi i)^n} \int_{B_1 \times \ldots \times
B_n} f(\lambda_1, \ldots, \lambda_n)
(\lambda_1 I - \Ra_1)^{-1} \ldots (\lambda_n I - \Ra_n)^{-1}
d\lambda_1 \ldots d\lambda_n$$
(compare to theorem VII.1.10 in \cite{DS}).

\end{document}